\documentclass{class}

\usepackage[colorlinks,urlcolor=blue,linkcolor=blue,citecolor=blue]{hyperref}

\usepackage{color,array}

\usepackage{graphicx}

\usepackage{xcolor}
\usepackage{graphicx}
\usepackage{amsmath,amssymb}
\usepackage{xcolor}
\usepackage{setspace} 
\usepackage{units}
\usepackage[absolute,overlay]{textpos}
\usepackage{multirow}
\usepackage{multicol}
\usepackage{amsthm}
\usepackage{bm}
\usepackage{subfigure}
\usepackage{calligra}
\newtheorem{assum}{Assumption}
\newtheorem{definition}{Definition}
\newtheorem{remark}{Remark}
\newtheorem{theorem}{Theorem}
\newtheorem{lemma}{Lemma}

\jvol{00}
\jnum{XX}
\paper{1234567}
\pubyear{2026}
\receiveddate{XX September 2026}
\accepteddate{XX October 2026}
\publisheddate{XX November 2026}
\currentdate{XX November 2026}
\doiinfo{Doi Number}

\begin{document}

\sptitle{Article Category}

\title{\!\!\!\!\!\!\!\!\!\!\!\!\!\mbox{Nonholonomic Event-Triggered Source Seeking}}

\editor{This paper was recommended by Associate Editor X. Y. Z.}

\author{VICTOR HUGO PEREIRA RODRIGUES\affilmark{1}, 
TIAGO ROUX OLIVEIRA\affilmark{1},  
MIROSLAV KRSTI{\' C}\affilmark{2}}  

\affil{Department of Electronics and Telecommunication Engineering, State University of Rio de Janeiro, Rio de Janeiro -- RJ, Brazil \\ (e-mail: victor.rodrigues@uerj.br and tiagoroux@uerj.br).} 
\affil{Department of Mechanical and Aerospace Engineering at the Jacobs School of Engineering, University of California San Diego, La Jolla -- CA, USA (e-mail: mkrstic@ucsd.edu).} 

\corresp{CORRESPONDING AUTHOR: V.~H.~P. Rodrigues (e-mail: \href{mailto:victor.rodrigues@uerj.br}{victor.rodrigues@uerj.br})}
\authornote{This study was financed in part by the Coordena{\c c}{\~a}o de Aperfei{\c c}oamento de Pessoal de N{\'i}vel Superior – Brasil (CAPES) – Finance Code 001. The authors also acknowledge the Brazilian Funding Agencies Conselho Nacional de Desenvolvimento Cient{\'i}fico e Tecnol{\'o}gico (CNPq) and Funda{\c c}{\~a}o de Amparo {\`a} Pesquisa do Estado do Rio de Janeiro (FAPERJ).}


\markboth{NONHOLONOMIC EVENT-TRIGGERED SOURCE SEEKING}{V.~H.~P. RODRIGUES {\itshape ET AL}.}

\begin{abstract}
This paper introduces an event-triggered source seeking control (ET-SSC)
strategy for autonomous vehicles modeled as nonholonomic unicycles.
The key feature of the proposed design is the separation between
the continuous probing motion required for gradient extraction and the
event-triggered update of the source-seeking correction. This separation
enables aperiodic computation and transmission of the correction while
the periodic probing signals remain continuously active. The triggering
mechanism acts on a low-pass-filtered gradient estimate, obtained by
high-pass (washout) filtering, demodulation, and low-pass filtering of
the measured signal, and combines a relative and an absolute threshold.
Stability is proved with a time-periodic Lyapunov function constructed
along the exponentially stable periodic solution of the filtered closed
loop, establishing practical exponential convergence of the vehicle to a
neighborhood of the unknown source. The guaranteed minimum dwell time is
strictly positive and independent of the probing frequency, which
excludes Zeno behavior. Unlike available Lie-bracket analyses of the forward-velocity design, the proposed averaging-based
analysis accommodates anisotropic quadratic source fields with distinct
spatial curvatures. The resulting ET-SSC framework therefore provides a
resource-aware implementation in which the source-seeking correction is
updated aperiodically while the information-generating probing motion
remains continuous.
\end{abstract}


\maketitle

\section{INTRODUCTION}

Networked control systems offer significant advantages in terms of installation and maintenance costs \cite{ZHGDDYP:2020}, but the limited bandwidth of the available communication channels may cause traffic congestion, transmission delays and packet dropouts \cite{HNX:2007}. Event-triggered controllers are an effective solution to this issue.

Event-Triggered (ET) control performs the control task on a non-periodic basis, responding to a triggering condition that is defined as a function of the plant's state. In addition to ensuring asymptotic stability \cite{T:2007}, this approach minimizes control effort by updating the control signal and transmitting data only when the error between the current state and the equilibrium set surpasses a threshold that could lead to instability \cite{BH:2013}. Many studies have extended event-based control to networked systems with high levels of complexity, addressing both linear \cite{HJT:2012} and nonlinear systems \cite{T:2007,APDN:2016}. More recent developments have also explored event-triggering concepts in nonlinear state estimation and learning-based control, where communication, computation, or data acquisition is activated only when required to preserve desired closed-loop properties \cite{NEW1,NEW}. From an information-theoretic perspective,\linebreak ET control has also been investigated by exploiting the timing of transmissions itself as a communication resource \cite{NEW2}. 

Despite the extensive body of literature on ET, as well as the established results of extremum seeking (ES) for static and general nonlinear dynamic systems in continuous time \cite{KW:2000}, the theoretical developments of ES have extended well beyond these systems. These advancements now encompass discrete-time \cite{CKAL:2002}, stochastic \cite{MK:2009, LK:2010} and multivariable systems \cite{GKN:2012}, noncooperative games \cite{FKB:2012}, time delays \cite{OKT:2017, ZFO:2023}, and systems governed by partial differential equations \cite{OK:2022_book_golden}. However, until recently, no work directly addressed the integration of ES with ET---this gap was filled by \cite{AUT:2025}, which developed multivariable ES algorithms based on perturbation (averaging) estimates of the gradient under ET control. The problem considered here lies at the intersection of this event-triggered ES framework and \textit{source seeking control} \cite{SCL_2007,source1,Azuma2012,Atanasov2015,source2,source3,source4,source5,source6}, where an autonomous vehicle must locate the source of an unknown scalar field in a GPS-denied environment using only local signal measurements.

Source seeking control (SSC) for nonholonomic vehicles has evolved along two main directions. The forward-velocity-tuning approach introduced in \cite{SCL_2007} employs a constant angular velocity together with a periodically perturbed linear velocity, producing a forward-backward probing motion whose averaged drift directs the vehicle toward the source. A second direction tunes the angular velocity while maintaining forward motion, beginning with \cite{source1} and followed by bounded-update-rate and Lie-bracket-based designs \cite{Scheinker2014,Durr2017,Grushkovskaya2018}. From the standpoint of vehicle motion, steering-based designs have important advantages: they avoid the repeated forward-backward motion of the original scheme and are applicable to vehicles, such as aircraft, that cannot reverse their direction of travel. Moreover, Lie-bracket averaging provides a particularly powerful framework for their analysis, often revealing the desired gradient flow directly from the interaction of the fast vector fields \cite{Durr2013,Durr2015,Durr2017}.

These advantages, however, expose a fundamental distinction when event-triggered implementation is considered. In Lie-bracket-based source seekers, the gradient information is generated intrinsically by the fast vehicle motion through the interaction of state-dependent oscillatory vector fields; there is generally no separately computed controller-side correction available to be held between events. Indeed, if the relevant state-dependent quantity is frozen over an inter-event interval, the corresponding vector fields become constant on that interval and their Lie bracket vanishes. Consequently, preserving the gradient-generating mechanism requires updates on the fast probing time scale. In contrast, the forward-velocity design considered here separates the continuous probing motion from the event-triggered updating of the controller-side correction. The periodic probing signals therefore remain uninterrupted while the correction signal is recomputed and transmitted only at triggering instants.

There is also an analytical reason for retaining classical periodic averaging in the present setting. The recent second-order Lie-bracket treatment of the forward-velocity design in \cite{source6} provides a modern interpretation of this class of seekers, but its analysis relies on distance-dependent source maps, for which the local Hessian is isotropic. Here, in contrast, the quadratic source field is allowed to be anisotropic, with distinct curvatures along the two spatial directions. For this more general setting, the rotating-coordinate invariance underlying the available second-order Lie-bracket analysis no longer holds, whereas classical averaging accommodates the anisotropic map directly.

Building on this observation, we extend the multivariable event-triggered ES framework of \cite{AUT:2025} to nonholonomic source seeking, with the vehicle and forward-velocity probing scheme of \cite{SCL_2007} and a vector gradient-demodulation structure related to multivariable ES \cite{GKN:2012}. The event-triggering mechanism acts specifically on the correction signal, so that the resource saving considered here concerns the computation and communication of the source-seeking correction, rather than a reduction of the vehicle's probing motion itself.

The proposed ET-SSC approach permits aperiodic correction updates. The gradient estimate is generated as in the classical extremum-seeking loop \cite{KW:2000} --- the measured signal is passed through a high-pass (washout) filter and demodulated --- and is then low-pass filtered; the triggering mechanism acts on this low-pass-filtered estimate, combining a relative threshold with a small absolute one. Because the filtered estimate does not inherit the fast rotation of the demodulation signal, the guaranteed minimum dwell time is strictly positive and independent of the probing frequency, thereby excluding the \emph{Zeno phenomenon}. Stability is proved with a time-periodic Lyapunov function constructed along the exponentially stable periodic solution of the filtered continuous-time closed loop, which yields practical exponential convergence of the vehicle position to a neighborhood of the unknown source. To the best of our knowledge, this provides the first event-triggered source-seeking framework for a nonholonomic vehicle in which the controller-side correction can be updated aperiodically with a guaranteed minimum inter-update time while the fast probing motion remains continuous.

\emph{Notation:} $\|\cdot\|$ denotes the Euclidean norm and its induced matrix norm, and $\lambda_{\min}(\cdot)$ and $\lambda_{\max}(\cdot)$ the minimum and maximum eigenvalues of a symmetric matrix. $\mathcal{O}(\cdot)$ denotes the standard big-$O$ notation: $\delta_{1}(\varepsilon)=\mathcal{O}(\delta_{2}(\varepsilon))$ if there exist positive constants $k$ and $c$ such that $|\delta_{1}(\varepsilon)| \le k |\delta_{2}(\varepsilon)|$ for all $|\varepsilon|<c$. For a symmetric matrix $X$, $X  > 0$ $(< 0)$ means that $X$ is positive (negative) definite.

\section{PROBLEM FORMULATION}

\subsection{Nonholonomic Mobile Robots}

\begin{figure}[h!]
	\centering
    \IfFileExists{blockDiagram_sourceSeeking.pdf}{\includegraphics[width=8.3cm]{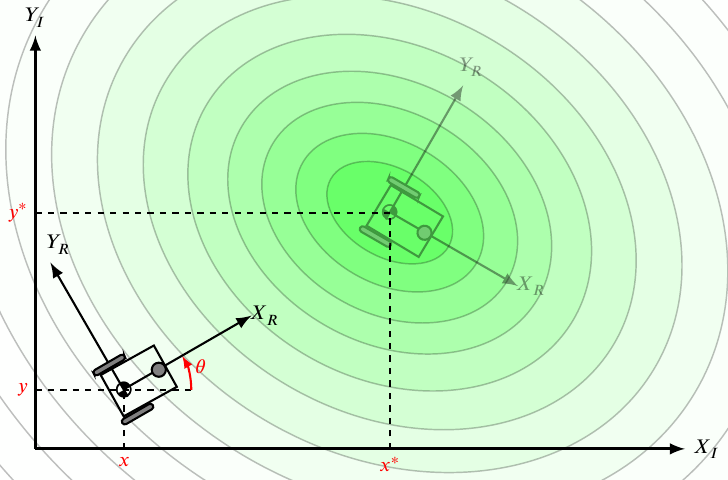}}%
    {\fbox{\parbox[c][4cm][c]{8cm}{\centering Figure file \texttt{blockDiagram\_sourceSeeking.pdf} missing}}}
    %
	\caption{Source-seeking problem for a nonholonomic mobile robot with position $(x,y)$ and heading angle $\theta$. The unknown source is located at $(x^*,y^*)$; the shades of green represent the source signal field.}
	\label{fig:unicycle}
\end{figure}

We consider a unicycle-type mobile robot equipped with a sensor, which may be placed directly at the vehicle's center or at some distance $r$ from it \cite{SCL_2007}. Fig.~\ref{fig:unicycle} illustrates the position and orientation of the robot.
The robot is modeled as a rigid body moving on a plane, with a global (inertial) frame $\{X_I, Y_I\}$ and a local frame $\{X_R, Y_R\}$ attached to the robot. Its pose is described by the position $(x(t),y(t))$ of its center in the global frame and by the heading angle $\theta(t)$ between the axes $X_R$ and $X_I$.

The motion of the robot center follows the equations
\begin{align}
    \frac{dx}{dt}(t) &= v(t) \cos (\theta(t)), \label{eq:dotX_ETSSC_v1} \\
    \frac{dy}{dt}(t) &= v(t) \sin (\theta(t)), \label{eq:dotY_ETSSC_v1} \\
    \frac{d\theta}{dt}(t) &= \omega(t), \label{eq:dotTheta_ETSSC_v1}
\end{align}
where the control inputs are the linear speed $v(t)$, directed along $X_R$, and the angular speed $\omega(t)$. The nonholonomic constraint prevents lateral sliding along $Y_R$.

\subsection{Source Signal Modeling}

The vehicle's objective is to locate a source emitting a signal whose strength decreases with distance from the source. We assume this signal field is described by an unknown nonlinear function
\begin{align}
f(x, y) &= f^{\ast} - \frac{f_{x}}{2}(x(t) - x^{\ast})^2 - \frac{f_{y}}{2}(y(t) - y^{\ast})^2 , \label{eq:Q_ETSSC_v1}
\end{align}
which admits a single local maximum $f^{\ast} = f(x^{\ast}, y^{\ast})$ at the point $(x^{\ast}, y^{\ast})$, representing the source location. The objective is to steer the vehicle toward this maximum using only the signal values measured at the sensor, under the following assumption.

\begin{assum}\label{assum:f}
The unique optimizer vector $[x^{\ast}\,, y^{\ast}]^{\top} \in \mathbb{R}^{2}$ and the scalar $f^{\ast} > 0$ are unknown parameters of the nonlinear map (\ref{eq:Q_ETSSC_v1}).
\end{assum}

\subsection{Linear Velocity $v(t)$ and Angular Velocity $\omega(t)$}

With the sensor placed at the robot's center, the linear and angular velocities are designed as
\begin{align}
    v(t)&=a_{1}\omega_{1}\cos(\omega_{1}t) +2L^{\top}(t)u(t)\,, \label{eq:vt_ETSSC_v1} \\
     \omega(t)&=\omega_{0} \label{eq:wt_ETSSC_v1}\,,
\end{align}
where 
\begin{align}
    L(t)&:=\begin{bmatrix}
        \cos(\omega_{0}t+\theta_{0}) \\
        \sin(\omega_{0}t+\theta_{0})
    \end{bmatrix} \,, \label{eq:L}
\end{align}
and $u(t) \in \mathbb{R}^{2}$ is ET-SSC law to be designed, $a_{1} >0$ is the amplitude of a small dither signal and $\omega_{0},\omega_{1}>0$ are the probing frequencies while $\theta_{0} \in \mathbb{R}$ is the known initial heading angle.

\begin{assum}\label{assum:w}
The probing frequencies satisfy
\begin{align}
\omega_1 = \mu  \omega_0\,, \label{eq:w1_w2}
\end{align}
for any positive integer $\mu>3$ \cite{SCL_2007}.
\end{assum} 

Substituting (\ref{eq:vt_ETSSC_v1}) into (\ref{eq:dotX_ETSSC_v1}) and (\ref{eq:dotY_ETSSC_v1}), and (\ref{eq:wt_ETSSC_v1}) into (\ref{eq:dotTheta_ETSSC_v1}), and using Assumption~\ref{assum:w}, the system dynamics (\ref{eq:dotX_ETSSC_v1})--(\ref{eq:dotTheta_ETSSC_v1}) can be rewritten as
\begin{align}
    \frac{dx}{dt}(t)&= \frac{a_{1}\mu \omega_{0}}{2}\cos((\mu +1)\omega_{0}t+\theta_{0}) \nonumber \\
    &\quad
    +\frac{a_{1}\mu \omega_{0}}{2}\cos((\mu -1)\omega_{0}t-\theta_{0}) \nonumber \\
    &\quad+2\cos(\omega_{0}t+\theta_{0})L^{\top}(t)u(t) \,, \label{eq:dotX_ETSSC_v2} \\
    \frac{dy}{dt}(t)&= \frac{a_{1}\mu \omega_{0}}{2}\sin((\mu +1)\omega_{0}t+\theta_{0}) \nonumber \\
    &\quad
    +\frac{a_{1}\mu \omega_{0}}{2}\sin((\mu -1)\omega_{0}t-\theta_{0}) \nonumber \\
    &\quad+2\sin(\omega_{0}t+\theta_{0})L^{\top}(t)u(t)\,, \label{eq:dotY_ETSSC_v2} \\
     \frac{d\theta}{dt}(t)&=\omega_{0} \Rightarrow \theta(t)=\omega_{0}t+\theta_{0} \,.\label{eq:dotTheta_ETSSC_v2}
\end{align}

\subsection{Estimation Error}

The first two terms on the right-hand sides of (\ref{eq:dotX_ETSSC_v2}) and (\ref{eq:dotY_ETSSC_v2}) are the time derivatives of the components of the dither vector
\begin{align}
S(t) &:= \left[S_{1}(t),\ S_{2}(t)\right]^\top, \label{eq:S_v1}\\
S_{1}(t)&=\frac{a_{1}}{2}\frac{\mu}{\mu+1}\sin((\mu+1)\omega_{0}t+\theta_{0}) \nonumber \\
    &\quad+\frac{a_{1}}{2}\frac{\mu}{\mu-1}\sin((\mu-1)\omega_{0}t-\theta_{0})\,, \label{eq:S1_v1} \\
S_{2}(t)&=-\frac{a_{1}}{2}\frac{\mu}{\mu+1}\cos((\mu+1)\omega_{0}t+\theta_{0}) \nonumber \\
    &\quad+\frac{a_{1}}{2}\frac{\mu}{\mu-1}\cos((\mu-1)\omega_{0}t-\theta_{0})\,.  \label{eq:S2_v1}
\end{align}
Hence, defining the position vector, the source vector, the position estimate and the estimation error as
\begin{align}
    p(t) &:= [x(t),\ y(t)]^\top \,, \quad p^\ast := [x^\ast,\ y^\ast]^\top  \,,  \label{eq:q_ETSSC_v1} \\
    \hat{p}(t) &:= p(t)-S(t) \,, \quad \tilde{p}(t) :=  \hat{p}(t)-p^{\ast} \,, \label{eq:tildeQ_ETSSC_v1}
\end{align}
one has
\begin{align}
p(t) - p^\ast = \tilde{p}(t) + S(t)\,. \label{eq:q_ETSSC_v2}
\end{align}  
Unlike \cite{source6}, we now define the \textit{anisotropic Hessian matrix} ($f_x \neq f_y$) of the nonlinear map (\ref{eq:Q_ETSSC_v1}) as
\begin{align}
H^{\ast} := \begin{bmatrix} -f_{x} & 0 \\ 0 & -f_{y} \end{bmatrix}, \label{eq:H1}
\end{align}
where $H^{\ast}=H^{\ast \top}<0$, such that, using (\ref{eq:q_ETSSC_v1}), the map (\ref{eq:Q_ETSSC_v1}) can be rewritten as
\begin{align}
f(p) = f^{\ast} + \frac{1}{2}(p(t) - p^{\ast})^{\top} H^{\ast} (p(t) - p^{\ast}). \label{eq:Q_ETSSC_v2}
\end{align}
By plugging (\ref{eq:q_ETSSC_v2}) into (\ref{eq:Q_ETSSC_v2}), the map $f$ can be rewritten with respect to the estimation error as
\begin{align}
f(t,\tilde{p}) &= f^{\ast} + \frac{1}{2}(\tilde{p}(t) +S(t))^{\top} H^{\ast} (\tilde{p}(t) +S(t)) \\
&=f^{\ast} + \frac{1}{2}\tilde{p}^{\top}(t)H^{\ast}\tilde{p}(t) +S^{\top}(t) H^{\ast} \tilde{p}(t) \nonumber \\
&\quad+ \frac{1}{2}S^{\top}(t)H^{\ast}S(t) \,. \label{eq:Q_ETSSC_v3}
\end{align}

On the other hand, from (\ref{eq:tildeQ_ETSSC_v1}), (\ref{eq:dotX_ETSSC_v2}), (\ref{eq:dotY_ETSSC_v2}) and (\ref{eq:L}), the dynamics that governs $\hat{p}(t)$, as well as $\tilde{p}(t)$, is given by
\begin{align}
\frac{d\tilde{p}}{dt}(t)&=\frac{d\hat{p}}{dt}(t)=2L(t)L^{\top}(t)u(t) \nonumber \\
&=(I_{2}+\Upsilon(t))u(t)\label{eq:gradient_dtildeThetadt_20250206_1}\,, 
\end{align}
where $I_{2}$ denotes the $2\times2$ identity matrix and $\Upsilon(t) \in \mathbb{R}^{2 \times 2}$ is given by
\begin{align}
    \Upsilon(t):=\begin{bmatrix}
        \cos(2\omega_{0}t+2\theta_{0}) & \sin(2\omega_{0}t+2\theta_{0})\\
        \sin(2\omega_{0}t+2\theta_{0}) & -\cos(2\omega_{0}t+2\theta_{0})
    \end{bmatrix} \,. \label{eq:Upsilon}
\end{align}
Since $2L(t)L^{\top}(t)=I_{2}+\Upsilon(t)$ and $\|L(t)\|=1$, the matrix $I_{2}+\Upsilon(t)$ has rank one and
\begin{align}
\|I_{2}+\Upsilon(t)\|=2\|L(t)\|^{2}=2\,, \quad \forall t\geq 0\,. \label{eq:F_normUps}
\end{align}

\subsection{Demodulation Vector and Dither Identities}

Throughout the paper, the signal measured by the vehicle is denoted by
\begin{align}
f(t):=f(p(t))\,. \label{eq:F_meas}
\end{align}
The gradient information contained in (\ref{eq:F_meas}) is extracted by means of the demodulation vector
\begin{align}
M(t)&=\left[M_{1}(t)\,,M_{2}(t)\right]^{\top} \,, \label{eq:M_ETSSC_v1}
\end{align}
with entries
\begin{align}
M_{1}(t)&= \frac{2(\mu+1)(\mu-1)^2}{a_{1}\mu(\mu^2+1)}\sin((\mu+1)\omega_{0}t+\theta_{0}) \nonumber \\
&\quad+ \frac{2(\mu+1)^2(\mu-1)}{a_{1}\mu(\mu^2+1)}\sin((\mu-1)\omega_{0}t-\theta_{0})\,, \label{eq:M1_ETSSC_v1} \\
M_{2}(t)&= -\frac{2(\mu+1)(\mu-1)^2}{a_{1}\mu(\mu^2+1)}\cos((\mu+1)\omega_{0}t+\theta_{0}) \nonumber \\
&\quad+ \frac{2(\mu+1)^2(\mu-1)}{a_{1}\mu(\mu^2+1)}\cos((\mu-1)\omega_{0}t-\theta_{0})\,. \label{eq:M2_ETSSC_v1}
\end{align}
To make the oscillatory nature of the term $M(t)S^{\top}(t)$ explicit, we use the demodulation vector components (\ref{eq:M_ETSSC_v1})--(\ref{eq:M2_ETSSC_v1}) to define the time-varying matrix $\Delta(t) \in \mathbb{R}^{2 \times 2}$ as
\begin{align}
    \Delta(t):=\begin{bmatrix}
        \Delta_{11}(t) & \Delta_{12}(t)\\
        \Delta_{21}(t) & \Delta_{22}(t)
    \end{bmatrix} \,, \label{eq:Delta}
\end{align}
whose entries are given by

\begin{align}
\Delta_{11}(t) &=- \frac{(\mu-1)^2}{2(\mu^2+1)}\cos\big(2(\mu+1)\omega_{0}t+2\theta_{0}\big) \nonumber \\
&\quad- \frac{(\mu+1)^2}{2(\mu^2+1)}\cos\big(2(\mu-1)\omega_{0}t-2\theta_{0}\big) \nonumber \\
&\quad + \frac{\mu^2-1}{\mu^2+1}\Big[\cos\big(2\omega_{0}t+2\theta_{0}\big)-\cos\big(2\mu\omega_{0}t\big)\Big]\,, \label{eq:Delta11}\\
\Delta_{12}(t) &= -\frac{(\mu-1)^2}{2(\mu^2+1)}\sin\big(2(\mu+1)\omega_{0}t+2\theta_{0}\big) \nonumber \\
&\quad+ \frac{(\mu+1)^2}{2(\mu^2+1)}\sin\big(2(\mu-1)\omega_{0}t-2\theta_{0}\big) \nonumber\\
&\quad + \frac{\mu^2-1}{\mu^2+1}\sin\big(2\omega_{0}t+2\theta_{0}\big)\,, \label{eq:Delta12}
\end{align}
\begin{align}
\Delta_{21}(t) &= \Delta_{12}(t)\,, \label{eq:Delta21}\\
\Delta_{22}(t) &=\frac{(\mu-1)^2}{2(\mu^2+1)}\cos\big(2(\mu+1)\omega_{0}t+2\theta_{0}\big) \nonumber \\
&\quad+ \frac{(\mu+1)^2}{2(\mu^2+1)}\cos\big(2(\mu-1)\omega_{0}t-2\theta_{0}\big) \nonumber\\
&\quad - \frac{\mu^2-1}{\mu^2+1}\Big[\cos\big(2\omega_{0}t+2\theta_{0}\big)+\cos\big(2\mu\omega_{0}t\big)\Big]\,, \label{eq:Delta22}
\end{align}
which allows us to write $M(t)S^{\top}(t)$ compactly as
\begin{align}
M(t)S^{\top}(t)& = I_{2}+\Delta(t)\,. \label{eq:MST}
\end{align}

All the signals introduced so far are generated by sinusoids oscillating at the two commensurate frequencies $\omega_{0}$ and $\omega_{1}=\mu\omega_{0}$ (with $\mu$ a positive integer). Hence, $L(t)$, $\Upsilon(t)$, $S(t)$, $M(t)$ and $\Delta(t)$ are periodic with the common period
\begin{align}
T&= 2\pi \times \text{LCM}\left\{\frac{1}{\omega_{0}}\,,\frac{1}{\mu\omega_{0}}\right\}= \frac{2\pi}{\omega_{0}}\,, \label{eq:T}
\end{align}
where LCM denotes the least common multiple. Moreover, they depend on $t$ and $\omega_{0}$ only through the product $\omega_{0}t$, so that their amplitudes are independent of $\omega_{0}$. For a continuous $T$-periodic function $h(t)$, scalar-, vector- or matrix-valued, its mean value is denoted by
\begin{align}
\overline{h}:=\frac{1}{T}\int_{0}^{T}h(\tau)d\tau\,. \label{eq:F_mean}
\end{align}

Finally, we isolate the term of (\ref{eq:Q_ETSSC_v3}) that does not depend on the estimation error $\tilde{p}(t)$, together with its mean value,
\begin{align}
s(t):=\frac{1}{2}S^{\top}(t)H^{\ast}S(t)\,, \quad s_{0}:=\overline{s}=\frac{1}{T}\int_{0}^{T}s(\tau)d\tau\,. \label{eq:F_s}
\end{align}
To compute them, let
\begin{align}
\phi_{\pm}(t):=(\mu\pm1)\omega_{0}t\pm\theta_{0}\,, \quad c_{\pm}:=\frac{a_{1}}{2}\frac{\mu}{\mu\pm1}\,, \label{eq:F_phi}
\end{align}
so that (\ref{eq:S1_v1})--(\ref{eq:S2_v1}) read $S_{1}=c_{+}\sin\phi_{+}+c_{-}\sin\phi_{-}$ and $S_{2}=-c_{+}\cos\phi_{+}+c_{-}\cos\phi_{-}$, with $\phi_{+}(t)-\phi_{-}(t)=2\omega_{0}t+2\theta_{0}$ and $\phi_{+}(t)+\phi_{-}(t)=2\mu\omega_{0}t$. Using $\sin^{2}\phi=\frac{1}{2}(1-\cos 2\phi)$, $\cos^{2}\phi=\frac{1}{2}(1+\cos 2\phi)$, $2\sin\phi_{+}\sin\phi_{-}=\cos(\phi_{+}-\phi_{-})-\cos(\phi_{+}+\phi_{-})$ and $2\cos\phi_{+}\cos\phi_{-}=\cos(\phi_{+}-\phi_{-})+\cos(\phi_{+}+\phi_{-})$, one obtains
\begin{align}
S_{1}^{2}(t)&=\frac{c_{+}^{2}+c_{-}^{2}}{2}-g(t)-c_{+}c_{-}\cos(2\mu\omega_{0}t)\,, \label{eq:F_S1sq}\\
S_{2}^{2}(t)&=\frac{c_{+}^{2}+c_{-}^{2}}{2}+g(t)-c_{+}c_{-}\cos(2\mu\omega_{0}t)\,, \label{eq:F_S2sq}\\
g(t)&:=\frac{c_{+}^{2}}{2}\cos\big(2(\mu+1)\omega_{0}t+2\theta_{0}\big) \nonumber \\
&\quad+\frac{c_{-}^{2}}{2}\cos\big(2(\mu-1)\omega_{0}t-2\theta_{0}\big) \nonumber \\
&\quad-c_{+}c_{-}\cos\big(2\omega_{0}t+2\theta_{0}\big)\,. \label{eq:F_g}
\end{align}
Since $s(t)=-\frac{1}{2}\big(f_{x}S_{1}^{2}(t)+f_{y}S_{2}^{2}(t)\big)$ by (\ref{eq:H1}), and since every cosine in (\ref{eq:F_S1sq})--(\ref{eq:F_g}) has zero mean over the period $T$, it follows that
\begin{align}
s_{0}&=-\frac{(f_{x}+f_{y})(c_{+}^{2}+c_{-}^{2})}{4} \nonumber \\
&=-\frac{(f_{x}+f_{y})a_{1}^{2}\mu^{2}(\mu^{2}+1)}{8(\mu^{2}-1)^{2}}\,, \label{eq:F_s0}\\
s(t)-s_{0}&=\frac{f_{x}+f_{y}}{2}c_{+}c_{-}\cos(2\mu\omega_{0}t)+\frac{f_{x}-f_{y}}{2}g(t)\,. \label{eq:F_ss0}
\end{align}
Therefore, with the quadratic form
\begin{align}
q(\tilde{p}):=\frac{1}{2}\tilde{p}^{\top}H^{\ast}\tilde{p}\,, \label{eq:F_q}
\end{align}
the measured signal (\ref{eq:Q_ETSSC_v3}) can be written as
\begin{align}
f(t)=f^{\ast}+s_{0}+q(\tilde{p}(t))+S^{\top}(t)H^{\ast}\tilde{p}(t)+\big(s(t)-s_{0}\big)\,. \label{eq:F_fdecomp}
\end{align}
Notice from (\ref{eq:F_ss0}) and (\ref{eq:F_g}) that $s(t)-s_{0}$ contains only the frequencies $2\omega_{0}$, $2(\mu-1)\omega_{0}$, $2\mu\omega_{0}$ and $2(\mu+1)\omega_{0}$, the first, second and fourth of which are present only in the anisotropic case $f_{x}\neq f_{y}$.

The next lemma collects uniform bounds on the periodic signals, which are used repeatedly in the sequel.

\begin{lemma}[Bounds on the dither signals] \label{lem:dither}
Under Assumption~\ref{assum:w}, for all $t\geq 0$,
\begin{align}
\|S(t)\|&\leq a_{\mu}:=\frac{\mu^{2}}{\mu^{2}-1}a_{1}\,, \label{eq:F_Sbound}\\
m\leq\|M(t)\|&\leq \mu\, m\,, \quad m:=\frac{4(\mu^{2}-1)}{a_{1}\mu(\mu^{2}+1)}\,, \label{eq:F_Mbound}\\
\|I_{2}+\Delta(t)\|&\leq \frac{4\mu^{2}}{\mu^{2}+1}<4\,, \label{eq:F_Dbound}\\
|s(t)-s_{0}|&\leq \|H^{\ast}\|a_{\mu}^{2}\,, \label{eq:F_sbound}
\end{align}
where $\|H^{\ast}\|=\max\{f_{x},f_{y}\}$.
\end{lemma}

\begin{proof}
With the notation (\ref{eq:F_phi}), $S(t)=c_{+}[\sin\phi_{+}\,,-\cos\phi_{+}]^{\top}+c_{-}[\sin\phi_{-}\,,\cos\phi_{-}]^{\top}$, where both vectors in brackets have unit norm and their inner product is $\sin\phi_{+}\sin\phi_{-}-\cos\phi_{+}\cos\phi_{-}=-\cos(\phi_{+}+\phi_{-})$. Hence,
\begin{align}
\|S(t)\|^{2}=c_{+}^{2}+c_{-}^{2}-2c_{+}c_{-}\cos(2\mu\omega_{0}t)\leq(c_{+}+c_{-})^{2}\,, \label{eq:F_Snorm}
\end{align}
where $c_{+}+c_{-}=\frac{a_{1}\mu}{2}\left(\frac{1}{\mu+1}+\frac{1}{\mu-1}\right)=\frac{\mu^{2}}{\mu^{2}-1}a_{1}$, which proves (\ref{eq:F_Sbound}). Analogously, (\ref{eq:M1_ETSSC_v1})--(\ref{eq:M2_ETSSC_v1}) read $M(t)=m_{+}[\sin\phi_{+}\,,-\cos\phi_{+}]^{\top}+m_{-}[\sin\phi_{-}\,,\cos\phi_{-}]^{\top}$, with
\begin{align}
m_{+}:=\frac{2(\mu+1)(\mu-1)^{2}}{a_{1}\mu(\mu^{2}+1)}\,, \quad m_{-}:=\frac{2(\mu+1)^{2}(\mu-1)}{a_{1}\mu(\mu^{2}+1)}\,, \label{eq:F_mpm}
\end{align}
so that $\|M(t)\|^{2}=m_{+}^{2}+m_{-}^{2}-2m_{+}m_{-}\cos(2\mu\omega_{0}t)$ and $(m_{-}-m_{+})^{2}\leq\|M(t)\|^{2}\leq(m_{+}+m_{-})^{2}$. Since $m_{-}-m_{+}=\frac{2(\mu^{2}-1)}{a_{1}\mu(\mu^{2}+1)}\big[(\mu+1)-(\mu-1)\big]=m$ and $m_{+}+m_{-}=\frac{2(\mu^{2}-1)}{a_{1}\mu(\mu^{2}+1)}\big[(\mu-1)+(\mu+1)\big]=\mu\, m$, inequality (\ref{eq:F_Mbound}) follows, with $m>0$ because $\mu>3$. By (\ref{eq:MST}), $I_{2}+\Delta(t)$ is the rank-one matrix $M(t)S^{\top}(t)$, whose induced norm is $\|M(t)\|\|S(t)\|\leq\mu\, m\, a_{\mu}=\frac{4\mu^{2}}{\mu^{2}+1}$, which is (\ref{eq:F_Dbound}). Finally, from (\ref{eq:F_s}) and (\ref{eq:F_Sbound}), $|s(t)|\leq\frac{1}{2}\|H^{\ast}\|a_{\mu}^{2}$ for all $t$ and, consequently, $|s_{0}|\leq\frac{1}{2}\|H^{\ast}\|a_{\mu}^{2}$, so that (\ref{eq:F_sbound}) follows from the triangle inequality.
\end{proof}

\section{CONTINUOUS-TIME GRADIENT-BASED TUNING LAW} \label{sec:CT}

\begin{figure*}[t!]
\centering
\includegraphics[width=14cm]{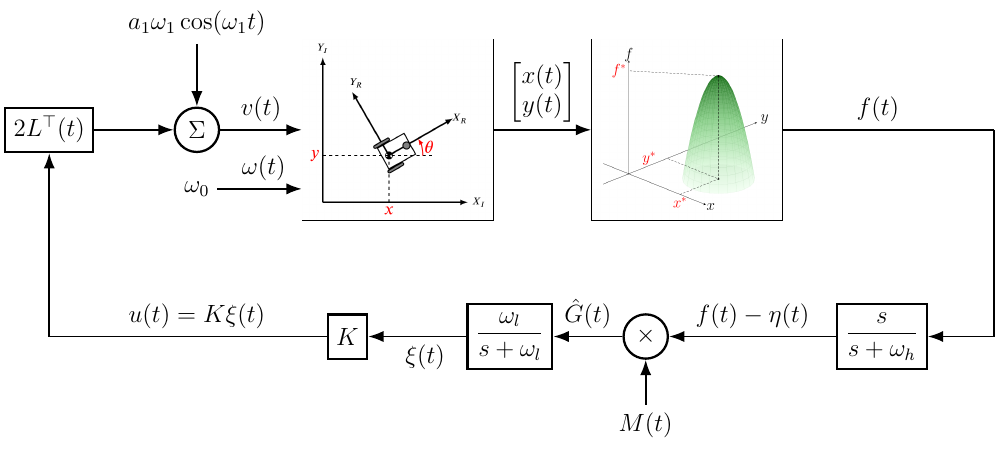}
\caption{Source Seeking Control for Nonholonomic Systems with high-pass and low-pass filters: continuous-time gradient-based tuning law.}
\label{fig:BD_SSC}
\end{figure*}

\subsection{High-Pass and Low-Pass Filtered Gradient Estimate}

As in the classical extremum-seeking loop \cite{KW:2000}, the measured signal (\ref{eq:F_meas}) is processed by a washout (high-pass) filter $s/(s+\omega_{h})$ before the demodulation stage and by a low-pass filter $\omega_{l}/(s+\omega_{l})$ after it, with cut-off frequencies $\omega_{h}>0$ and $\omega_{l}>0$. Let $\eta(t)\in\mathbb{R}$ be the state of the high-pass filter, governed by
\begin{align}
\frac{d\eta}{dt}(t)=-\omega_{h}\eta(t)+\omega_{h}f(t)\,. \label{eq:F_eta}
\end{align}
Since (\ref{eq:F_eta}) is a realization of the transfer function $\omega_{h}/(s+\omega_{h})$ from $f$ to $\eta$, and $1-\omega_{h}/(s+\omega_{h})=s/(s+\omega_{h})$, the output of the high-pass filter is $f(t)-\eta(t)$. This signal is multiplied by the demodulation vector $M(t)$, which yields the gradient estimate
\begin{align}
\hat{G}(t):=M(t)\big(f(t)-\eta(t)\big)\,. \label{eq:hatG_ETSSC_v1}
\end{align}
The gradient estimate is then applied to the low-pass filter, whose output $\xi(t)\in\mathbb{R}^{2}$ satisfies
\begin{align}
\frac{d\xi}{dt}(t)=-\omega_{l}\xi(t)+\omega_{l}\hat{G}(t)\,. \label{eq:F_xi}
\end{align}
The vector $\xi(t)$ is the low-pass-filtered gradient estimate, and the continuous-time gradient-based tuning law, shown in the block diagram of Fig.~\ref{fig:BD_SSC}, is
\begin{align}
u(t)=K\xi(t)\,, \quad K=\text{diag}(k_{1},k_{2})>0\,. \label{eq:u}
\end{align}

Two features of the filters are essential in what follows. First, the washout filter removes from the measurement (\ref{eq:F_fdecomp}) its constant component $f^{\ast}+s_{0}$, which would otherwise appear in $\hat{G}(t)$ multiplied by $M(t)$, as an oscillation that does not decay as the vehicle approaches the source. Second, the right-hand side of (\ref{eq:F_xi}) contains $M(t)$ but not its time derivative, so that the rate of change of $\xi(t)$ is bounded independently of $\omega_{0}$; this is what allows, in Section~\ref{sec:ETSSC}, a minimum dwell time that does not shrink as $\omega_{0}$ grows.

\subsection{Closed-Loop System}

Define, for analysis only, the shifted filter state and the closed-loop state vector
\begin{align}
\tilde{\eta}(t)&:=\eta(t)-f^{\ast}-s_{0}\,, \label{eq:F_etatilde}\\
z(t)&:=\big[\tilde{p}^{\top}(t)\,,\ \xi^{\top}(t)\,,\ \tilde{\eta}(t)\big]^{\top}\in\mathbb{R}^{5}\,. \label{eq:F_z}
\end{align}
Subtracting $\eta(t)=\tilde{\eta}(t)+f^{\ast}+s_{0}$ from (\ref{eq:F_fdecomp}), the output of the high-pass filter is
\begin{align}
f(t)-\eta(t)&=q(\tilde{p}(t))+S^{\top}(t)H^{\ast}\tilde{p}(t) \nonumber \\
&\quad+\big(s(t)-s_{0}\big)-\tilde{\eta}(t)\,. \label{eq:F_hp}
\end{align}
Multiplying (\ref{eq:F_hp}) by $M(t)$ and using the identity (\ref{eq:MST}) in the term $M(t)S^{\top}(t)H^{\ast}\tilde{p}(t)$, the gradient estimate (\ref{eq:hatG_ETSSC_v1}) becomes
\begin{align}
\hat{G}(t)&=M(t)q(\tilde{p}(t))+\left(H^{\ast}+\Delta(t)H^{\ast}\right)\tilde{p}(t) \nonumber \\
&\quad+M(t)\big(s(t)-s_{0}\big)-M(t)\tilde{\eta}(t)\,. \label{eq:hatG_ETSSC_v3}
\end{align}
The second term on the right-hand side of (\ref{eq:hatG_ETSSC_v3}) contains $H^{\ast}\tilde{p}(t)$, which is the gradient of the map (\ref{eq:Q_ETSSC_v2}) evaluated at $p^{\ast}+\tilde{p}(t)$; the remaining terms are either quadratic in $\tilde{p}(t)$, or have zero mean, or are multiplied by the zero-mean vector $M(t)$, as detailed in Section~\ref{sec:average}.

The three components of $z(t)$ evolve as follows. First, plugging (\ref{eq:u}) into (\ref{eq:gradient_dtildeThetadt_20250206_1}),
\begin{align}
\frac{d\tilde{p}}{dt}(t)=(I_{2}+\Upsilon(t))K\xi(t)\,. \label{eq:F_dp}
\end{align}
Second, substituting (\ref{eq:hatG_ETSSC_v3}) into (\ref{eq:F_xi}),
\begin{align}
\frac{d\xi}{dt}(t)&=\omega_{l}(I_{2}+\Delta(t))H^{\ast}\tilde{p}(t)-\omega_{l}\xi(t)-\omega_{l}M(t)\tilde{\eta}(t) \nonumber \\
&\quad+\omega_{l}M(t)\big(s(t)-s_{0}\big)+\omega_{l}M(t)q(\tilde{p}(t))\,. \label{eq:F_dxi}
\end{align}
Third, since $f^{\ast}$ and $s_{0}$ are constants, $\frac{d\tilde{\eta}}{dt}=\frac{d\eta}{dt}=\omega_{h}(f-\eta)$ by (\ref{eq:F_eta}), and (\ref{eq:F_hp}) gives
\begin{align}
\frac{d\tilde{\eta}}{dt}(t)&=\omega_{h}S^{\top}(t)H^{\ast}\tilde{p}(t)-\omega_{h}\tilde{\eta}(t) \nonumber \\
&\quad+\omega_{h}\big(s(t)-s_{0}\big)+\omega_{h}q(\tilde{p}(t))\,. \label{eq:F_deta}
\end{align}
Collecting (\ref{eq:F_dp})--(\ref{eq:F_deta}), the closed-loop system is, exactly,
\begin{align}
\frac{dz}{dt}(t)=A(t)z(t)+b(t)+N(t,z(t))\,, \label{eq:F_cl}
\end{align}
where
\begin{align}
A(t)&:=\begin{bmatrix}
0_{2\times 2} & (I_{2}+\Upsilon(t))K & 0_{2\times 1}\\
\omega_{l}(I_{2}+\Delta(t))H^{\ast} & -\omega_{l}I_{2} & -\omega_{l}M(t)\\
\omega_{h}S^{\top}(t)H^{\ast} & 0_{1\times 2} & -\omega_{h}
\end{bmatrix}, \label{eq:F_A}\\
n(t)&:=\begin{bmatrix}
0_{2\times 1}\\ \omega_{l}M(t)\\ \omega_{h}
\end{bmatrix}, \label{eq:F_n}\\
b(t)&:=n(t)\big(s(t)-s_{0}\big)\,, \quad N(t,z):=n(t)\,q(\tilde{p})\,. \label{eq:F_bN}
\end{align}
The matrix $A(t)$ and the vectors $n(t)$ and $b(t)$ are $T$-periodic and depend on $t$ and $\omega_{0}$ only through $\omega_{0}t$; hence their amplitudes are independent of $\omega_{0}$. Since $b(t)\not\equiv 0$, $z=0$ is not an equilibrium of (\ref{eq:F_cl}): the steady state is the periodic solution constructed in Lemma~\ref{lem:periodic}.

For later use, we bound the terms of (\ref{eq:F_cl}). The induced norm of a block matrix does not exceed the sum of the norms of its blocks; thus, from (\ref{eq:F_normUps}) and Lemma~\ref{lem:dither}, for all $t$,
\begin{align}
\|A(t)\|&\leq L_{A}:=2\|K\|+\omega_{l}\big(4\|H^{\ast}\|+1+\mu\, m\big) \nonumber \\
&\qquad\qquad+\omega_{h}\big(a_{\mu}\|H^{\ast}\|+1\big)\,, \label{eq:F_LA}\\
\|n(t)\|&\leq\omega_{l}\,\mu\, m+\omega_{h}\,. \label{eq:F_nbound}
\end{align}
Moreover, $|q(\tilde{p})|\leq\frac{1}{2}\|H^{\ast}\|\|\tilde{p}\|^{2}$ and, because $H^{\ast}$ is symmetric, $q(\tilde{p}_{a})-q(\tilde{p}_{b})=\frac{1}{2}(\tilde{p}_{a}+\tilde{p}_{b})^{\top}H^{\ast}(\tilde{p}_{a}-\tilde{p}_{b})$. Therefore, with
\begin{align}
C_{N}:=\frac{1}{2}\big(\omega_{l}\,\mu\, m+\omega_{h}\big)\|H^{\ast}\|\,, \label{eq:F_CN}
\end{align}
and recalling (\ref{eq:F_sbound}) and that $\|\tilde{p}\|\leq\|z\|$, one has, for all $t$ and all $z,z_{a},z_{b}\in\mathbb{R}^{5}$,
\begin{align}
\|b(t)\|&\leq 2C_{N}a_{\mu}^{2}\,, \label{eq:F_bbound}\\
\|N(t,z)\|&\leq C_{N}\|z\|^{2}\,, \label{eq:F_Nbound}\\
\|N(t,z_{a})-N(t,z_{b})\|&\leq C_{N}\big(\|z_{a}\|+\|z_{b}\|\big)\|z_{a}-z_{b}\|\,. \label{eq:F_Nlip}
\end{align}

\subsection{Average System} \label{sec:average}

We now compute the mean values (\ref{eq:F_mean}) of the periodic coefficients of (\ref{eq:F_cl}).
\begin{enumerate}
\item[(i)] $\overline{\Upsilon}=0$, since every entry of (\ref{eq:Upsilon}) is a sinusoid of frequency $2\omega_{0}$.
\item[(ii)] $\overline{S}=0$ and $\overline{M}=0$, since the entries of (\ref{eq:S1_v1})--(\ref{eq:S2_v1}) and (\ref{eq:M1_ETSSC_v1})--(\ref{eq:M2_ETSSC_v1}) are sums of sinusoids of frequencies $(\mu+1)\omega_{0}$ and $(\mu-1)\omega_{0}$, both nonzero.
\item[(iii)] $\overline{\Delta}=0$, since the entries (\ref{eq:Delta11})--(\ref{eq:Delta22}) are sums of sinusoids of frequencies $2(\mu+1)\omega_{0}$, $2(\mu-1)\omega_{0}$, $2\omega_{0}$ and $2\mu\omega_{0}$, all nonzero.
\item[(iv)] $\overline{s-s_{0}}=0$, by the definition of $s_{0}$ in (\ref{eq:F_s}).
\item[(v)] $\overline{M(s-s_{0})}=0$. Indeed, the product of two sinusoids of frequencies $n_{1}\omega_{0}$ and $n_{2}\omega_{0}$ is a sum of two sinusoids of frequencies $(n_{1}+n_{2})\omega_{0}$ and $(n_{1}-n_{2})\omega_{0}$, and thus has zero mean whenever $n_{1}\neq n_{2}$. The entries of $M(t)$ contain the frequencies $n_{1}\in\{\mu-1,\mu+1\}$, whereas, by (\ref{eq:F_ss0}) and (\ref{eq:F_g}), $s(t)-s_{0}$ contains the frequencies $n_{2}\in\{2,2\mu-2,2\mu,2\mu+2\}$. For a positive integer $\mu$, the equalities $\mu+1=2$, $\mu+1=2\mu$ and $\mu-1=2\mu-2$ hold only for $\mu=1$; the equalities $\mu+1=2\mu-2$ and $\mu-1=2$ hold only for $\mu=3$; and $\mu+1=2\mu+2$, $\mu-1=2\mu$, $\mu-1=2\mu+2$ are impossible. Hence, the two sets of frequencies are disjoint for every integer $\mu>3$. This is the point where Assumption~\ref{assum:w} is needed: for $\mu=3$ and $f_{x}\neq f_{y}$ the mean of $M(s-s_{0})$ does not vanish.
\end{enumerate}
From (i)--(v) and (\ref{eq:F_A})--(\ref{eq:F_bN}), $\overline{b}=\overline{n(s-s_{0})}=0$ and
\begin{align}
\bar{A}:=\overline{A}=\begin{bmatrix}
0_{2\times 2} & K & 0_{2\times 1}\\
\omega_{l}H^{\ast} & -\omega_{l}I_{2} & 0_{2\times 1}\\
0_{1\times 2} & 0_{1\times 2} & -\omega_{h}
\end{bmatrix}, \quad \overline{n}=\begin{bmatrix}
0_{2\times 1}\\ 0_{2\times 1}\\ \omega_{h}
\end{bmatrix}. \label{eq:F_Abar}
\end{align}
Accordingly, the average system associated with (\ref{eq:F_cl}) is
\begin{align}
\frac{dz_{\rm av}}{dt}(t)=\bar{A}z_{\rm av}(t)+\overline{n}\,q(\tilde{p}_{\rm av}(t))\,, \label{eq:F_av}
\end{align}
with $z_{\rm av}=[\tilde{p}_{\rm av}^{\top},\xi_{\rm av}^{\top},\tilde{\eta}_{\rm av}]^{\top}$. Notice that the Hessian $H^{\ast}$ enters the average system directly, with no restriction on $f_{x}$ and $f_{y}$ other than their positivity.

Writing $f_{1}:=f_{x}$ and $f_{2}:=f_{y}$, and reordering the state as $(\tilde{p}_{1},\xi_{1})$, $(\tilde{p}_{2},\xi_{2})$, $\tilde{\eta}$, the matrix $\bar{A}$ becomes block diagonal with the blocks
\begin{align}
\bar{A}_{i}:=\begin{bmatrix}
0 & k_{i}\\ -\omega_{l}f_{i} & -\omega_{l}
\end{bmatrix}, \quad i=1,2\,, \label{eq:F_Ai}
\end{align}
and the scalar $-\omega_{h}$. The eigenvalues of $\bar{A}$ are therefore $-\omega_{h}$ and the roots of
\begin{align}
\lambda^{2}+\omega_{l}\lambda+\omega_{l}k_{i}f_{i}=0\,, \quad i=1,2\,, \label{eq:F_charpoly}
\end{align}
which have negative real parts because the coefficients $\omega_{l}$ and $\omega_{l}k_{i}f_{i}$ are positive. Thus, $\bar{A}$ is Hurwitz for every diagonal $K>0$ and every $\omega_{l},\omega_{h}>0$, without any additional assumption on the gains. The solution $\bar{P}=\bar{P}^{\top}$ of the Lyapunov equation
\begin{align}
\bar{A}^{\top}\bar{P}+\bar{P}\bar{A}=-I_{5} \label{eq:F_Pbar}
\end{align}
is available in closed form: in the original ordering (\ref{eq:F_z}) of the state,
\begin{align}
\bar{P}&=\begin{bmatrix}
\text{diag}(\bar{p}_{11,1},\bar{p}_{11,2}) & \text{diag}(\bar{p}_{12,1},\bar{p}_{12,2}) & 0_{2\times 1}\\
\text{diag}(\bar{p}_{12,1},\bar{p}_{12,2}) & \text{diag}(\bar{p}_{22,1},\bar{p}_{22,2}) & 0_{2\times 1}\\
0_{1\times 2} & 0_{1\times 2} & \frac{1}{2\omega_{h}}
\end{bmatrix}, \label{eq:F_Pbar_explicit}\\
\bar{p}_{12,i}&:=\frac{1}{2\omega_{l}f_{i}}\,, \quad \bar{p}_{22,i}:=\frac{\omega_{l}f_{i}+k_{i}}{2\omega_{l}^{2}f_{i}}\,, \nonumber \\
\bar{p}_{11,i}&:=\frac{1}{k_{i}}\left(\frac{\omega_{l}f_{i}+k_{i}}{2\omega_{l}}+\frac{1}{2f_{i}}\right). \label{eq:F_pij}
\end{align}
Indeed, for each block (\ref{eq:F_Ai}), with $\bar{P}_{i}:=\left[\begin{smallmatrix}\bar{p}_{11,i} & \bar{p}_{12,i}\\ \bar{p}_{12,i} & \bar{p}_{22,i}\end{smallmatrix}\right]$, the entries $(1,1)$, $(2,2)$ and $(1,2)$ of the symmetric matrix $\bar{A}_{i}^{\top}\bar{P}_{i}+\bar{P}_{i}\bar{A}_{i}$ are, respectively, $-2\omega_{l}f_{i}\,\bar{p}_{12,i}=-1$, $2k_{i}\,\bar{p}_{12,i}-2\omega_{l}\,\bar{p}_{22,i}=-1$, $k_{i}\,\bar{p}_{11,i}-\omega_{l}\,\bar{p}_{12,i}-\omega_{l}f_{i}\,\bar{p}_{22,i}=0$, while, for the scalar block, $2(-\omega_{h})\frac{1}{2\omega_{h}}=-1$. Moreover, $\bar{P}>0$, since $\bar{p}_{11,i}>0$ and $\bar{p}_{11,i}\,\bar{p}_{22,i}-\bar{p}_{12,i}^{2}=\frac{(\omega_{l}f_{i}+k_{i})^{2}+\omega_{l}^{2}}{4\omega_{l}^{3}k_{i}f_{i}}>0$. We denote
\begin{align}
\gamma_{1}:=\lambda_{\min}(\bar{P})>0\,, \quad \gamma_{2}:=\lambda_{\max}(\bar{P})=\|\bar{P}\|\,. \label{eq:F_gamma}
\end{align}
Two elementary consequences of (\ref{eq:F_Pbar}) are used below. Since $\bar{A}$ is the mean of $A(t)$, $\|\bar{A}\|\leq\frac{1}{T}\int_{0}^{T}\|A(\tau)\|d\tau\leq L_{A}$; and taking norms in (\ref{eq:F_Pbar}), $1=\|\bar{A}^{\top}\bar{P}+\bar{P}\bar{A}\|\leq 2\|\bar{A}\|\|\bar{P}\|$, that is,
\begin{align}
\|\bar{A}\|\leq L_{A}\,, \quad 2L_{A}\gamma_{2}\geq 1\,. \label{eq:F_LAgamma}
\end{align}

\subsection{Auxiliary Results}

The relation between the closed-loop system (\ref{eq:F_cl}) and its average, which is the subject of the general averaging theory \cite{K:2002}, is established here directly, by means of four lemmas with explicit constants. Throughout, a constant is said to be \emph{independent of $\omega_{0}$} if it depends only on $\mu$, $a_{1}$, $f_{x}$, $f_{y}$, $K$, $\omega_{l}$ and $\omega_{h}$.

\begin{lemma}[Zero-mean antiderivative] \label{lem:antiderivative}
Let $h(t)$ be a continuous, $T$-periodic, matrix-valued function with zero mean, $\overline{h}=0$. Then
\begin{align}
\mathcal{I}[h](t):=\frac{1}{T}\int_{t-T}^{t}\left(\int_{\nu}^{t}h(\tau)d\tau\right)d\nu \label{eq:F_Idef}
\end{align}
is continuously differentiable, $T$-periodic, has zero mean, and satisfies, for all $t$,
\begin{align}
\frac{d}{dt}\mathcal{I}[h](t)&=h(t)\,, \label{eq:F_Ider}\\
\|\mathcal{I}[h](t)\|&\leq\frac{T}{2}\sup_{\tau}\|h(\tau)\|=\frac{\pi}{\omega_{0}}\sup_{\tau}\|h(\tau)\|\,. \label{eq:F_Ibound}
\end{align}
\end{lemma}

\begin{proof}
Let $h_{0}(t):=\int_{0}^{t}h(\tau)d\tau$, which is $T$-periodic because $\int_{t}^{t+T}h(\tau)d\tau=T\,\overline{h}=0$. Since $\int_{\nu}^{t}h(\tau)d\tau=h_{0}(t)-h_{0}(\nu)$ and the integral of the $T$-periodic function $h_{0}$ over any interval of length $T$ equals $T\,\overline{h_{0}}$, definition (\ref{eq:F_Idef}) gives $\mathcal{I}[h](t)=h_{0}(t)-\overline{h_{0}}$. Hence $\mathcal{I}[h]$ is continuously differentiable with derivative $h(t)$, is $T$-periodic, and has zero mean. Finally, for $\nu\in[t-T,t]$, $\left\|\int_{\nu}^{t}h(\tau)d\tau\right\|\leq(t-\nu)\sup_{\tau}\|h(\tau)\|$, so that $\|\mathcal{I}[h](t)\|\leq\frac{\sup_{\tau}\|h(\tau)\|}{T}\int_{t-T}^{t}(t-\nu)d\nu=\frac{T}{2}\sup_{\tau}\|h(\tau)\|$, and $T=2\pi/\omega_{0}$ by (\ref{eq:T}).
\end{proof}

\begin{lemma}[Uniform exponential stability of the linear part] \label{lem:transition}
Let $E(t)\in\mathbb{R}^{5\times 5}$ be continuous and satisfy $\sup_{t}\|E(t)\|\leq\frac{1}{24\gamma_{2}}$, and let $\Phi_{E}(t,\tau)$ denote the transition matrix of
\begin{align}
\frac{d\chi}{dt}(t)=\big(A(t)+E(t)\big)\chi(t)\,. \label{eq:F_lin}
\end{align}
If
\begin{align}
\omega_{0}\geq 128\pi\gamma_{2}L_{A}^{2}\,, \label{eq:F_w1}
\end{align}
then, for all $t\geq\tau$,
\begin{align}
\|\Phi_{E}(t,\tau)\|\leq k_{\Phi}e^{-\lambda_{\Phi}(t-\tau)}\,, ~~ k_{\Phi}:=3\sqrt{\frac{\gamma_{2}}{\gamma_{1}}}\,, ~~ \lambda_{\Phi}:=\frac{1}{4\gamma_{2}}\,, \label{eq:F_Phibound}
\end{align}
where $k_{\Phi}$ and $\lambda_{\Phi}$ are independent of $\omega_{0}$ and of $E(t)$.
\end{lemma}

\begin{proof}
\emph{Step 1: near-identity transformation.} By Section~\ref{sec:average}, $\tilde{A}(t):=A(t)-\bar{A}$ is continuous, $T$-periodic and has zero mean, and $\|\tilde{A}(t)\|\leq 2L_{A}$ by (\ref{eq:F_LA}) and (\ref{eq:F_LAgamma}). Let $W(t):=\mathcal{I}[\tilde{A}](t)$. By Lemma~\ref{lem:antiderivative},
\begin{align}
\frac{dW}{dt}(t)=\tilde{A}(t)\,, \quad \|W(t)\|\leq\frac{2\pi L_{A}}{\omega_{0}}\leq\frac{1}{2}\,, \label{eq:F_W}
\end{align}
where the last inequality holds because (\ref{eq:F_w1}) and (\ref{eq:F_LAgamma}) give $\omega_{0}\geq 64\pi L_{A}(2L_{A}\gamma_{2})\geq 64\pi L_{A}$. Consequently, $I_{5}+W(t)$ is invertible, with
\begin{align}
\|I_{5}+W(t)\|&\leq\frac{3}{2}\,, \nonumber \\
\|(I_{5}+W(t))^{-1}\|&\leq\frac{1}{1-\|W(t)\|}\leq 2\,. \label{eq:F_IW}
\end{align}

\emph{Step 2: transformed system.} Let $\chi(t)=(I_{5}+W(t))w(t)$. Differentiating and using (\ref{eq:F_W}), $\frac{d\chi}{dt}=\tilde{A}w+(I_{5}+W)\frac{dw}{dt}$, whereas (\ref{eq:F_lin}) gives $\frac{d\chi}{dt}=(\bar{A}+\tilde{A}+E)(I_{5}+W)w$. Equating both expressions, the term $\tilde{A}w$ cancels and $(I_{5}+W)\frac{dw}{dt}=\big(\bar{A}+\bar{A}W+\tilde{A}W+E(I_{5}+W)\big)w$. Multiplying by $(I_{5}+W)^{-1}$ and using $(I_{5}+W)^{-1}\bar{A}=\bar{A}-(I_{5}+W)^{-1}W\bar{A}$,
\begin{align}
\frac{dw}{dt}(t)&=\big(\bar{A}+R(t)+\tilde{E}(t)\big)w(t)\,, \label{eq:F_w}\\
R&:=(I_{5}+W)^{-1}\big(\bar{A}W-W\bar{A}+\tilde{A}W\big)\,, \nonumber \\
\tilde{E}&:=(I_{5}+W)^{-1}E(I_{5}+W)\,. \nonumber
\end{align}
From (\ref{eq:F_W}), (\ref{eq:F_IW}), (\ref{eq:F_LAgamma}) and (\ref{eq:F_w1}),
\begin{align}
\|R(t)\|&\leq 2\big(2L_{A}+2L_{A}\big)\frac{2\pi L_{A}}{\omega_{0}}=\frac{16\pi L_{A}^{2}}{\omega_{0}}\leq\frac{1}{8\gamma_{2}}\,, \label{eq:F_R}\\
\|\tilde{E}(t)\|&\leq 2\cdot\frac{3}{2}\cdot\frac{1}{24\gamma_{2}}=\frac{1}{8\gamma_{2}}\,. \label{eq:F_Etilde}
\end{align}

\emph{Step 3: Lyapunov estimate.} Let $\upsilon(t):=w^{\top}(t)\bar{P}w(t)$. By (\ref{eq:F_w}), (\ref{eq:F_Pbar}), (\ref{eq:F_R}) and (\ref{eq:F_Etilde}), $\frac{d\upsilon}{dt}=w^{\top}\big(\bar{A}^{\top}\bar{P}+\bar{P}\bar{A}\big)w+2w^{\top}\bar{P}\big(R+\tilde{E}\big)w \leq-\|w\|^{2}+2\gamma_{2}\left(\frac{1}{8\gamma_{2}}+\frac{1}{8\gamma_{2}}\right)\|w\|^{2}=-\frac{1}{2}\|w\|^{2}$. Since $\gamma_{1}\|w\|^{2}\leq\upsilon\leq\gamma_{2}\|w\|^{2}$, it follows that $\frac{d\upsilon}{dt}\leq-\frac{1}{2\gamma_{2}}\upsilon$, i.e., $\frac{d}{dt}\big(e^{t/(2\gamma_{2})}\upsilon(t)\big)\leq 0$, so that $\upsilon(t)\leq e^{-(t-\tau)/(2\gamma_{2})}\upsilon(\tau)$ and $\|w(t)\|\leq\sqrt{\frac{\gamma_{2}}{\gamma_{1}}}\,e^{-\frac{t-\tau}{4\gamma_{2}}}\|w(\tau)\|\,, \quad \forall t\geq\tau$.

\emph{Step 4: original coordinates.} By (\ref{eq:F_IW}), $\|\chi(t)\|\leq\frac{3}{2}\|w(t)\|$ and $\|w(\tau)\|\leq 2\|\chi(\tau)\|$. Hence $\|\chi(t)\|\leq 3\sqrt{\gamma_{2}/\gamma_{1}}\,e^{-(t-\tau)/(4\gamma_{2})}\|\chi(\tau)\|$ for every initial condition $\chi(\tau)$, which is (\ref{eq:F_Phibound}).
\end{proof}

In the sequel, $\Phi_{A}(t,\tau):=\Phi_{E}(t,\tau)$ with $E\equiv 0$ denotes the transition matrix of $\frac{d\chi}{dt}=A(t)\chi$. We recall that, since $A(t)$ is $T$-periodic, $\Phi_{A}(t+T,\tau+T)=\Phi_{A}(t,\tau)$, and that $\frac{\partial}{\partial t}\Phi_{A}(t,\tau)=A(t)\Phi_{A}(t,\tau)$ and $\frac{\partial}{\partial\tau}\Phi_{A}(t,\tau)=-\Phi_{A}(t,\tau)A(\tau)$, the latter being obtained by differentiating $\Phi_{A}(t,\tau)\Phi_{A}(\tau,t)=I_{5}$ with respect to $\tau$.

\begin{lemma}[Periodic solution] \label{lem:periodic}
Define the constants, independent of $\omega_{0}$,
\begin{align}
C_{p}:=\frac{8\pi}{3}C_{N}a_{\mu}^{2}\left(1+\frac{k_{\Phi}L_{A}}{\lambda_{\Phi}}\right), \quad \rho_{p}:=\frac{\lambda_{\Phi}}{4k_{\Phi}C_{N}}\,. \label{eq:F_Cp}
\end{align}
If $\omega_{0}\geq\max\{128\pi\gamma_{2}L_{A}^{2}\,,\,C_{p}/\rho_{p}\}$, then the closed-loop system (\ref{eq:F_cl}) has a $T$-periodic solution $z_{\rm per}(t)=[\tilde{p}_{\rm per}^{\top}(t),\xi_{\rm per}^{\top}(t),\tilde{\eta}_{\rm per}(t)]^{\top}$ such that
\begin{align}
\sup_{t}\|z_{\rm per}(t)\|\leq\frac{C_{p}}{\omega_{0}}\,, \label{eq:F_zper_bound}
\end{align}
and $z_{\rm per}$ is the only $T$-periodic solution of (\ref{eq:F_cl}) satisfying $\sup_{t}\|z(t)\|\leq\rho_{p}$.
\end{lemma}

\begin{proof}
Let $\mathcal{X}$ be the set of continuous $T$-periodic functions $\zeta:\mathbb{R}\to\mathbb{R}^{5}$ with $\|\zeta\|_{\infty}:=\sup_{t}\|\zeta(t)\|\leq\rho_{p}$, which is a closed subset of the Banach space of continuous $T$-periodic functions with the norm $\|\cdot\|_{\infty}$. For $\zeta\in\mathcal{X}$ define
\begin{align}
(\mathcal{T}\zeta)(t):=\int_{-\infty}^{t}\Phi_{A}(t,\tau)\big(b(\tau)+N(\tau,\zeta(\tau))\big)d\tau\,, \label{eq:F_Top}
\end{align}
where the integral converges absolutely by (\ref{eq:F_Phibound}), (\ref{eq:F_bbound}) and (\ref{eq:F_Nbound}).

\emph{Step 1: $\mathcal{T}\zeta$ is $T$-periodic.} Changing the integration variable to $\tau-T$ in $(\mathcal{T}\zeta)(t+T)$ and using $\Phi_{A}(t+T,\tau+T)=\Phi_{A}(t,\tau)$ together with the $T$-periodicity of $b$, $n$ and $\zeta$ yields $(\mathcal{T}\zeta)(t+T)=(\mathcal{T}\zeta)(t)$.

\emph{Step 2: fixed points of $\mathcal{T}$ are the $T$-periodic solutions in $\mathcal{X}$.} Differentiating (\ref{eq:F_Top}) with $\Phi_{A}(t,t)=I_{5}$ and $\frac{\partial}{\partial t}\Phi_{A}(t,\tau)=A(t)\Phi_{A}(t,\tau)$, $\frac{d}{dt}(\mathcal{T}\zeta)(t)=b(t)+N(t,\zeta(t))+A(t)(\mathcal{T}\zeta)(t)$, so every fixed point $\zeta=\mathcal{T}\zeta$ solves (\ref{eq:F_cl}). Conversely, if $\zeta\in\mathcal{X}$ solves (\ref{eq:F_cl}), the variation-of-constants formula gives, for $t_{a}\leq t$, $\zeta(t)=\Phi_{A}(t,t_{a})\zeta(t_{a})+\int_{t_{a}}^{t}\Phi_{A}(t,\tau)\big(b(\tau)+N(\tau,\zeta(\tau))\big)d\tau$; letting $t_{a}\to-\infty$, the first term vanishes because $\|\Phi_{A}(t,t_{a})\zeta(t_{a})\|\leq k_{\Phi}e^{-\lambda_{\Phi}(t-t_{a})}\rho_{p}$, and thus $\zeta=\mathcal{T}\zeta$.

\emph{Step 3: the forced response is of order $1/\omega_{0}$.} Since $\overline{b}=0$, Lemma~\ref{lem:antiderivative} applies to $\beta(t):=\mathcal{I}[b](t)$, and (\ref{eq:F_bbound}) gives $\frac{d\beta}{dt}=b$ and $\|\beta(t)\|\leq 2\pi C_{N}a_{\mu}^{2}/\omega_{0}$. Since $\frac{\partial}{\partial\tau}\big(\Phi_{A}(t,\tau)\beta(\tau)\big)=\Phi_{A}(t,\tau)b(\tau)-\Phi_{A}(t,\tau)A(\tau)\beta(\tau)$, integrating over $\tau\in(-\infty,t]$ and noting that $\Phi_{A}(t,\tau)\beta(\tau)\to 0$ as $\tau\to-\infty$, $\int_{-\infty}^{t}\Phi_{A}(t,\tau)b(\tau)d\tau=\beta(t)+\int_{-\infty}^{t}\Phi_{A}(t,\tau)A(\tau)\beta(\tau)d\tau$. Hence, by (\ref{eq:F_Phibound}), (\ref{eq:F_LA}) and (\ref{eq:F_Cp}),
\begin{align}
\left\|\int_{-\infty}^{t}\Phi_{A}(t,\tau)b(\tau)d\tau\right\|&\leq\frac{2\pi C_{N}a_{\mu}^{2}}{\omega_{0}}\left(1+\frac{k_{\Phi}L_{A}}{\lambda_{\Phi}}\right) \nonumber \\
&=\frac{3}{4}\frac{C_{p}}{\omega_{0}}\,. \label{eq:F_forced}
\end{align}

\emph{Step 4: $\mathcal{T}$ maps $\mathcal{X}$ into itself and is a contraction.} For $\zeta\in\mathcal{X}$, by (\ref{eq:F_forced}), (\ref{eq:F_Phibound}), (\ref{eq:F_Nbound}), the definition of $\rho_{p}$ and $C_{p}/\omega_{0}\leq\rho_{p}$, $\|\mathcal{T}\zeta\|_{\infty}\leq\frac{3}{4}\frac{C_{p}}{\omega_{0}}+\frac{k_{\Phi}}{\lambda_{\Phi}}C_{N}\rho_{p}^{2}=\frac{3}{4}\frac{C_{p}}{\omega_{0}}+\frac{\rho_{p}}{4}\leq\rho_{p}$. For $\zeta_{a},\zeta_{b}\in\mathcal{X}$, by (\ref{eq:F_Nlip}), $\|\mathcal{T}\zeta_{a}-\mathcal{T}\zeta_{b}\|_{\infty}\leq\frac{k_{\Phi}}{\lambda_{\Phi}}C_{N}\,2\rho_{p}\|\zeta_{a}-\zeta_{b}\|_{\infty} =\frac{1}{2}\|\zeta_{a}-\zeta_{b}\|_{\infty}$. Therefore, the successive approximations $\zeta_{j+1}=\mathcal{T}\zeta_{j}$, $\zeta_{0}=0$, converge uniformly to the unique fixed point $z_{\rm per}$ of $\mathcal{T}$ in $\mathcal{X}$ (contraction mapping principle \cite{K:2002}). By Steps~1 and~2, $z_{\rm per}$ is a $T$-periodic solution of (\ref{eq:F_cl}) and it is the only one in $\mathcal{X}$.

\emph{Step 5: size of $z_{\rm per}$.} From $z_{\rm per}=\mathcal{T}z_{\rm per}$, (\ref{eq:F_forced}), (\ref{eq:F_Nbound}) and $\|z_{\rm per}\|_{\infty}\leq\rho_{p}$, $\|z_{\rm per}\|_{\infty}\leq\frac{3}{4}\frac{C_{p}}{\omega_{0}}+\frac{k_{\Phi}}{\lambda_{\Phi}}C_{N}\rho_{p}\|z_{\rm per}\|_{\infty} =\frac{3}{4}\frac{C_{p}}{\omega_{0}}+\frac{1}{4}\|z_{\rm per}\|_{\infty}$, which gives (\ref{eq:F_zper_bound}).
\end{proof}

The last lemma is stated for the closed loop driven by an additional input, which will be the sampling-induced error in Section~\ref{sec:ETSSC}; the continuous-time loop (\ref{eq:F_cl}) corresponds to $e\equiv 0$. Let
\begin{align}
B(t):=\begin{bmatrix}
(I_{2}+\Upsilon(t))K\\ 0_{2\times 2}\\ 0_{1\times 2}
\end{bmatrix}, \quad \|B(t)\|\leq C_{B}:=2\|K\|\,, \label{eq:F_B}
\end{align}
where the bound follows from (\ref{eq:F_normUps}).

\begin{lemma}[Error dynamics and periodic Lyapunov matrix] \label{lem:lyapunov}
Let
\begin{align}
\omega_{0}\geq\max\big\{128\pi\gamma_{2}L_{A}^{2}\,,\ C_{p}/\rho_{p}\,,\ 48\gamma_{2}C_{N}C_{p}\big\}\,, \label{eq:F_w3}
\end{align}
and let $z(t)$ be a solution, on some interval, of
\begin{align}
\frac{dz}{dt}(t)=A(t)z(t)+b(t)+N(t,z(t))+B(t)e(t)\,, \label{eq:F_cl_e}
\end{align}
with $e(t)\in\mathbb{R}^{2}$ piecewise continuous. Then:
\begin{enumerate}
\item[(a)] The error $\tilde{z}(t):=z(t)-z_{\rm per}(t)$ satisfies, exactly,
\begin{align}
\frac{d\tilde{z}}{dt}(t)&=A_{\rm per}(t)\tilde{z}(t)+N(t,\tilde{z}(t))+B(t)e(t)\,, \label{eq:F_derr}\\
A_{\rm per}(t)&:=A(t)+E_{\rm per}(t)\,, \nonumber \\
E_{\rm per}(t)&:=n(t)\,\tilde{p}_{\rm per}^{\top}(t)H^{\ast}\big[I_{2}~~0_{2\times 3}\big]\,, \label{eq:F_Eper}
\end{align}
where $A_{\rm per}(t)$ is $T$-periodic and
\begin{align}
\|E_{\rm per}(t)\|\leq\frac{2C_{N}C_{p}}{\omega_{0}}\leq\frac{1}{24\gamma_{2}}\,. \label{eq:F_Eper_bound}
\end{align}
\item[(b)] Let $\Phi_{\rm per}(\tau,t)$ be the transition matrix of $\frac{d\chi}{dt}=A_{\rm per}(t)\chi$. The matrix
\begin{align}
P(t):=\int_{t}^{\infty}\Phi_{\rm per}^{\top}(\tau,t)\Phi_{\rm per}(\tau,t)d\tau \label{eq:F_P}
\end{align}
is symmetric, continuously differentiable, $T$-periodic, and satisfies, for all $t$,

\begin{align}
&\frac{dP}{dt}(t)+A_{\rm per}^{\top}(t)P(t)+P(t)A_{\rm per}(t)=-I_{5}\,, \label{eq:F_Pdot}\\
&\alpha_{1}I_{5}\leq P(t)\leq\alpha_{2}I_{5}\,, \label{eq:F_Pbounds}\\
&\alpha_{1}:=\frac{1}{2L_{A}+\frac{1}{12\gamma_{2}}}\,, \quad \alpha_{2}:=\frac{k_{\Phi}^{2}}{2\lambda_{\Phi}}=\frac{18\gamma_{2}^{2}}{\gamma_{1}}\,, \label{eq:F_alphas}
\end{align}
with $\alpha_{1}$ and $\alpha_{2}$ independent of $\omega_{0}$.
\item[(c)] At every $t$ where $e(t)$ is continuous, the function $V(t):=\tilde{z}^{\top}(t)P(t)\tilde{z}(t)$ satisfies
\begin{align}
\frac{dV}{dt}&\leq-\|\tilde{z}\|^{2}+2\alpha_{2}C_{N}\|\tilde{z}\|^{3}+2\alpha_{2}C_{B}\|\tilde{z}\|\|e\|\,. \label{eq:F_Vdot}
\end{align}
\end{enumerate}
\end{lemma}

\begin{proof}
\emph{(a)} Let $\tilde{z}_{p}$ denote the first two components of $\tilde{z}$, so that $\tilde{p}=\tilde{p}_{\rm per}+\tilde{z}_{p}$. Since $H^{\ast}$ is symmetric, $q(\tilde{p}_{\rm per}+\tilde{z}_{p})=q(\tilde{p}_{\rm per})+\tilde{p}_{\rm per}^{\top}H^{\ast}\tilde{z}_{p}+q(\tilde{z}_{p})$, and (\ref{eq:F_bN}) gives $N(t,z_{\rm per}+\tilde{z})=N(t,z_{\rm per})+n(t)\tilde{p}_{\rm per}^{\top}H^{\ast}\tilde{z}_{p}+N(t,\tilde{z})$. Subtracting $\frac{dz_{\rm per}}{dt}=A z_{\rm per}+b+N(t,z_{\rm per})$ from (\ref{eq:F_cl_e}), the term $b(t)$ cancels and (\ref{eq:F_derr}) follows, with $n(t)\tilde{p}_{\rm per}^{\top}H^{\ast}\tilde{z}_{p}=E_{\rm per}(t)\tilde{z}$. The matrix $A_{\rm per}(t)$ is $T$-periodic because $A(t)$, $n(t)$ and $z_{\rm per}(t)$ are. By (\ref{eq:F_nbound}), (\ref{eq:F_CN}) and (\ref{eq:F_zper_bound}), $\|E_{\rm per}(t)\|\leq\|n(t)\|\|H^{\ast}\|\|\tilde{p}_{\rm per}(t)\|\leq 2C_{N}C_{p}/\omega_{0}$, and the second inequality in (\ref{eq:F_Eper_bound}) is the third condition in (\ref{eq:F_w3}).

\emph{(b)} By (\ref{eq:F_Eper_bound}), Lemma~\ref{lem:transition} applies with $E=E_{\rm per}$, so that $\|\Phi_{\rm per}(\tau,t)\|\leq k_{\Phi}e^{-\lambda_{\Phi}(\tau-t)}$ for $\tau\geq t$. Hence the integral (\ref{eq:F_P}) converges and, for every $x\in\mathbb{R}^{5}$, $x^{\top}P(t)x=\int_{t}^{\infty}\|\Phi_{\rm per}(\tau,t)x\|^{2}d\tau\leq\frac{k_{\Phi}^{2}}{2\lambda_{\Phi}}\|x\|^{2}=\alpha_{2}\|x\|^{2}$. For the lower bound, let $\chi(\tau):=\Phi_{\rm per}(\tau,t)x$ and $L_{\rm per}:=L_{A}+\frac{1}{24\gamma_{2}}\geq\sup_{\tau}\|A_{\rm per}(\tau)\|$. Then $\frac{d}{d\tau}\|\chi\|^{2}=2\chi^{\top}A_{\rm per}\chi\geq-2L_{\rm per}\|\chi\|^{2}$, so that $\|\chi(\tau)\|^{2}\geq e^{-2L_{\rm per}(\tau-t)}\|x\|^{2}$ and $x^{\top}P(t)x\geq\int_{t}^{\infty}e^{-2L_{\rm per}(\tau-t)}d\tau\,\|x\|^{2}=\frac{\|x\|^{2}}{2L_{\rm per}}=\alpha_{1}\|x\|^{2}$. The $T$-periodicity of $P$ follows from $\Phi_{\rm per}(\tau+T,t+T)=\Phi_{\rm per}(\tau,t)$ after the change of variable $\tau\to\tau+T$ in (\ref{eq:F_P}). Finally, differentiating (\ref{eq:F_P}) with respect to $t$, using $\Phi_{\rm per}(t,t)=I_{5}$ and $\frac{\partial}{\partial t}\Phi_{\rm per}(\tau,t)=-\Phi_{\rm per}(\tau,t)A_{\rm per}(t)$,
\begin{align}
\frac{dP}{dt}(t)&=-I_{5}-A_{\rm per}^{\top}(t)\int_{t}^{\infty}\Phi_{\rm per}^{\top}(\tau,t)\Phi_{\rm per}(\tau,t)d\tau \nonumber \\
&\quad-\int_{t}^{\infty}\Phi_{\rm per}^{\top}(\tau,t)\Phi_{\rm per}(\tau,t)d\tau\,A_{\rm per}(t)\,,
\end{align}
which is (\ref{eq:F_Pdot}).

\emph{(c)} Along (\ref{eq:F_derr}), using (\ref{eq:F_Pdot}), $\frac{dV}{dt}=\tilde{z}^{\top}\Big(\frac{dP}{dt}+A_{\rm per}^{\top}P+PA_{\rm per}\Big)\tilde{z} +2\tilde{z}^{\top}P\big(N(t,\tilde{z})+Be\big) =-\|\tilde{z}\|^{2}+2\tilde{z}^{\top}PN(t,\tilde{z})+2\tilde{z}^{\top}PBe$, and (\ref{eq:F_Vdot}) follows from (\ref{eq:F_Pbounds}), (\ref{eq:F_Nbound}) and (\ref{eq:F_B}).
\end{proof}

\subsection{Stability of the Continuous-Time Loop}

The next theorem establishes the stability properties of the continuous-time closed loop. Let $\sigma_{0}\in(0,1)$ be an arbitrary constant. It is not a design parameter of the controller, but a degree of freedom of the analysis: it is the fraction of the decay rate of the Lyapunov function used to dominate the terms of the closed loop which are not quadratic in the error, and it sets a tradeoff between the guaranteed decay rate and the size of the guaranteed region of attraction. Define the constants, independent of $\omega_{0}$,
\begin{align}
\rho_{0}&:=\frac{\sigma_{0}}{4\alpha_{2}C_{N}}\,, \quad \rho:=\frac{\rho_{0}}{4}\sqrt{\frac{\alpha_{1}}{\alpha_{2}}}\,, \label{eq:F_dbar}\\
\omega_{0}^{\ast}&:=\max\Big\{128\pi\gamma_{2}L_{A}^{2}\,,\ \frac{C_{p}}{\rho_{p}}\,,\ 48\gamma_{2}C_{N}C_{p}\,,\ \frac{C_{p}}{\rho}\Big\}\,. \label{eq:F_w0star}
\end{align}

\begin{theorem} \label{thm:ssc}
Consider the closed-loop system (\ref{eq:F_cl}), composed of the vehicle (\ref{eq:dotX_ETSSC_v1})--(\ref{eq:dotTheta_ETSSC_v1}) with the velocities (\ref{eq:vt_ETSSC_v1})--(\ref{eq:wt_ETSSC_v1}), the filters (\ref{eq:F_eta})--(\ref{eq:F_xi}) and the continuous-time tuning law (\ref{eq:u}), under Assumptions~\ref{assum:f} and~\ref{assum:w}. Then:
\begin{enumerate}
\item[(a)] For every $\omega_{0}>0$, the origin $z_{\rm av}=0$ is a locally exponentially stable equilibrium of the average system (\ref{eq:F_av}).
\item[(b)] For every $\omega_{0}\geq\omega_{0}^{\ast}$, the closed-loop system (\ref{eq:F_cl}) has a locally exponentially stable $T$-periodic solution $z_{\rm per}(t)$, which satisfies $\sup_{t}\|z_{\rm per}(t)\|\leq C_{p}/\omega_{0}$ and is the unique $T$-periodic solution in the ball $\sup_{t}\|z(t)\|\leq\rho_{p}$.
\item[(c)] For every $\omega_{0}\geq\omega_{0}^{\ast}$ and every initial condition with $\|z(0)\|\leq\rho$, the solution exists for all $t\geq 0$ and
\begin{align}
\|p(t)-p^{\ast}\|&\leq c_{1}e^{-\lambda t}\|z(0)\| \nonumber \\
&\quad+\mathcal{O}\left(\frac{\mu^{2}}{\mu^{2}-1}a_{1}+\frac{1}{\omega_{0}}\right), \label{eq:norm_p_thm1}
\end{align}
where $c_{1}:=\sqrt{\alpha_{2}/\alpha_{1}}$, the decay rate is
\begin{align}
\lambda:=\frac{1}{2}\frac{1}{\alpha_{2}}\big(1-\sigma_{0}\big)\,, \label{eq:F_lambda_c}
\end{align}
and the constant implicit in $\mathcal{O}(\cdot)$ is independent of $\omega_{0}$.
\end{enumerate}
\end{theorem}

\begin{proof}
\emph{(a)} The average system (\ref{eq:F_av}) does not depend on $\omega_{0}$. Let $\upsilon_{\rm av}:=z_{\rm av}^{\top}\bar{P}z_{\rm av}$. By (\ref{eq:F_Pbar}), and since $\|\overline{n}\,q(\tilde{p}_{\rm av})\|\leq\frac{\omega_{h}}{2}\|H^{\ast}\|\|z_{\rm av}\|^{2}$, $\frac{d\upsilon_{\rm av}}{dt}=-\|z_{\rm av}\|^{2}+2z_{\rm av}^{\top}\bar{P}\,\overline{n}\,q(\tilde{p}_{\rm av}) \leq-\|z_{\rm av}\|^{2}+\gamma_{2}\omega_{h}\|H^{\ast}\|\|z_{\rm av}\|^{3}$. Thus $\frac{d\upsilon_{\rm av}}{dt}\leq-\frac{1}{2}\|z_{\rm av}\|^{2}\leq-\frac{1}{2\gamma_{2}}\upsilon_{\rm av}$ whenever $\|z_{\rm av}\|\leq\rho_{\rm av}:=\frac{1}{2\gamma_{2}\omega_{h}\|H^{\ast}\|}$. The sublevel set $\{\upsilon_{\rm av}\leq\gamma_{1}\rho_{\rm av}^{2}\}$ is contained in the ball $\|z_{\rm av}\|\leq\rho_{\rm av}$ and is therefore forward invariant, and on it $\|z_{\rm av}(t)\|\leq\sqrt{\gamma_{2}/\gamma_{1}}\,e^{-t/(4\gamma_{2})}\|z_{\rm av}(0)\|$.

\emph{(b)--(c)} Let $\omega_{0}\geq\omega_{0}^{\ast}$. The existence, the size and the uniqueness of $z_{\rm per}$ are given by Lemma~\ref{lem:periodic}. Let $z(t)$ be a solution of (\ref{eq:F_cl}) defined on its maximal interval $[0,t_{\max})$, let $\tilde{z}(t)=z(t)-z_{\rm per}(t)$ and $V(t)=\tilde{z}^{\top}(t)P(t)\tilde{z}(t)$ as in Lemma~\ref{lem:lyapunov}, and suppose that
\begin{align}
\|\tilde{z}(0)\|\leq 2\rho=\frac{\rho_{0}}{2}\sqrt{\frac{\alpha_{1}}{\alpha_{2}}}\,. \label{eq:F_d0}
\end{align}
By (\ref{eq:F_Vdot}) with $e\equiv 0$, $\frac{dV}{dt}\leq-\|\tilde{z}\|^{2}+2\alpha_{2}C_{N}\|\tilde{z}\|^{3}$. Whenever $\|\tilde{z}(t)\|\leq\rho_{0}$, the definition of $\rho_{0}$ in (\ref{eq:F_dbar}) gives $2\alpha_{2}C_{N}\|\tilde{z}\|^{3}\leq\frac{\sigma_{0}}{2}\|\tilde{z}\|^{2}\leq\sigma_{0}\|\tilde{z}\|^{2}$, and therefore, by (\ref{eq:F_Pbounds}) and (\ref{eq:F_lambda_c}),
\begin{align}
\frac{dV}{dt}\leq-\big(1-\sigma_{0}\big)\|\tilde{z}\|^{2}\leq-\frac{1-\sigma_{0}}{\alpha_{2}}V=-2\lambda V\,. \label{eq:F_Vdot_c}
\end{align}
Let $t_{1}:=\sup\{t\in[0,t_{\max}):\ \|\tilde{z}(\tau)\|<\rho_{0}~~\forall\tau\in[0,t]\}$, which is positive because $\|\tilde{z}(0)\|\leq\rho_{0}/2$. On $[0,t_{1})$, (\ref{eq:F_Vdot_c}) shows that $V$ is nonincreasing, so that, by (\ref{eq:F_Pbounds}) and (\ref{eq:F_d0}), $\alpha_{1}\|\tilde{z}(t)\|^{2}\leq V(t)\leq V(0)\leq\alpha_{2}\|\tilde{z}(0)\|^{2}\leq\alpha_{1}\rho_{0}^{2}/4$, i.e., $\|\tilde{z}(t)\|\leq\rho_{0}/2$. If $t_{1}<t_{\max}$, continuity would give $\|\tilde{z}(t_{1})\|\leq\rho_{0}/2<\rho_{0}$ and $\|\tilde{z}\|$ would remain below $\rho_{0}$ on a right neighborhood of $t_{1}$, contradicting the definition of $t_{1}$. Hence $t_{1}=t_{\max}$; since $z(t)$ is then bounded on $[0,t_{\max})$, it can be continued and $t_{\max}=\infty$. Integrating $\frac{d}{dt}\big(e^{2\lambda t}V(t)\big)\leq 0$ and using (\ref{eq:F_Pbounds}),
\begin{align}
\|\tilde{z}(t)\|\leq\sqrt{\frac{\alpha_{2}}{\alpha_{1}}}\,e^{-\lambda t}\|\tilde{z}(0)\|\,, \quad \forall t\geq 0\,, \label{eq:F_d_c}
\end{align}
which proves that $z_{\rm per}$ is locally exponentially stable and completes the proof of (b).

For (c), $\|z(0)\|\leq\rho$ and $\|z_{\rm per}(0)\|\leq C_{p}/\omega_{0}\leq\rho$, by (\ref{eq:F_zper_bound}) and the last term in (\ref{eq:F_w0star}), imply (\ref{eq:F_d0}). From (\ref{eq:q_ETSSC_v2}), $\|\tilde{p}\|\leq\|z\|\leq\|\tilde{z}\|+\|z_{\rm per}\|$, (\ref{eq:F_Sbound}), (\ref{eq:F_zper_bound}) and (\ref{eq:F_d_c}),
\begin{align}
\|p(t)-p^{\ast}\|&\leq\|\tilde{p}(t)\|+\|S(t)\| \nonumber \\
&\leq\|\tilde{z}(t)\|+\frac{C_{p}}{\omega_{0}}+a_{\mu} \nonumber \\
&\leq c_{1}e^{-\lambda t}\Big(\|z(0)\|+\frac{C_{p}}{\omega_{0}}\Big)+\frac{C_{p}}{\omega_{0}}+a_{\mu} \nonumber \\
&\leq c_{1}e^{-\lambda t}\|z(0)\|+\frac{\mu^{2}}{\mu^{2}-1}a_{1}+\frac{(c_{1}+1)C_{p}}{\omega_{0}}\,, \label{eq:F_p_c}
\end{align}
where $e^{-\lambda t}\leq 1$ and $a_{\mu}=\frac{\mu^{2}}{\mu^{2}-1}a_{1}$ were used in the last step. Since $(c_{1}+1)C_{p}$ is independent of $\omega_{0}$, the last two terms of (\ref{eq:F_p_c}) are bounded by $\max\{1,(c_{1}+1)C_{p}\}\big(\frac{\mu^{2}}{\mu^{2}-1}a_{1}+\frac{1}{\omega_{0}}\big)$, which is (\ref{eq:norm_p_thm1}).
\end{proof}

\begin{remark}[Initialization and measured signal] \label{rem:init}
The states $\xi$ and $\eta$ are controller states and, by (\ref{eq:q_ETSSC_v2}), (\ref{eq:F_etatilde}) and Lemma~\ref{lem:dither},
\begin{align}
\|z(0)\|\leq\|p(0)-p^{\ast}\|+a_{\mu}+\|\xi(0)\|+|\eta(0)-f^{\ast}-s_{0}|\,. \label{eq:F_z0gen}
\end{align}
With $\xi(0)=0$ and $\eta(0)=f(p(0))$, the last term equals $|\frac{1}{2}(p(0)-p^{\ast})^{\top}H^{\ast}(p(0)-p^{\ast})-s_{0}|\leq\frac{\|H^{\ast}\|}{2}\big(\|p(0)-p^{\ast}\|^{2}+a_{\mu}^{2}\big)$, so that $\|z(0)\|$ depends only on the initial distance to the source and on the dither amplitude; with $\eta(0)=0$ it equals $|f^{\ast}+s_{0}|$, the transient of the washout filter, which decays with the time constant $1/\omega_{h}$. Finally, since $|f(t)-f^{\ast}|\leq\frac{\|H^{\ast}\|}{2}\|p(t)-p^{\ast}\|^{2}$ by (\ref{eq:Q_ETSSC_v2}), the bound (\ref{eq:norm_p_thm1}) also quantifies the convergence of the measured signal to a neighborhood of $f^{\ast}$ of order $\mathcal{O}\big(\frac{\mu^{4}}{(\mu^{2}-1)^{2}}a_{1}^{2}+\frac{1}{\omega_{0}^{2}}\big)$.
\end{remark}

\section{EVENT-TRIGGERED SOURCE SEEKING CONTROL} \label{sec:ETSSC}

\begin{figure*}[h!]
\centering
\includegraphics[width=15cm]{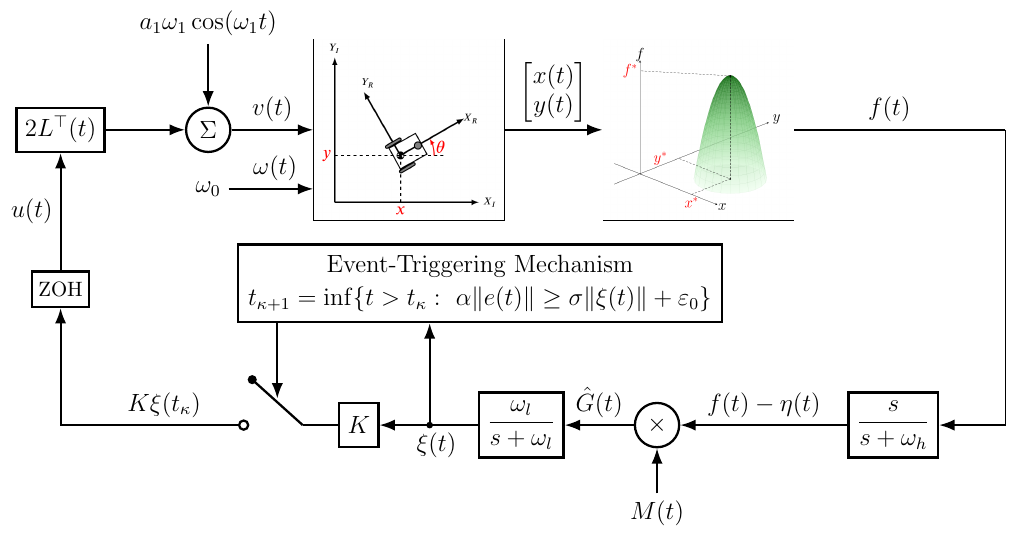}
\caption{Event-Triggered Source Seeking Control for Nonholonomic Systems with high-pass and low-pass filters: event-triggered tuning law.}
\label{fig:BD_NewtonES_v5}
\end{figure*}

We now propose a static event-triggered implementation that emulates the continuous-time tuning law (\ref{eq:u}) while reducing the frequency of control updates. As illustrated in Fig.~\ref{fig:BD_NewtonES_v5}, which differs from the continuous-time scheme of Fig.~\ref{fig:BD_SSC} only by the sampler, the zero-order hold and the event-triggering mechanism placed between the low-pass filter and the vehicle, the key idea is to monitor a suitably designed, error-dependent triggering condition that determines when the control law must be updated, rather than enforcing updates at every time instant. Control updates are executed only at a sequence of time instants $(t_{\kappa})_{\kappa \in \mathbb{N}}$, generated online by the event-triggering mechanism, which is designed so that the closed-loop system preserves the stability properties of its continuous-time counterpart.

\subsection{Emulation of the Continuous-Time Tuning Law}

Assuming no delays in the sensor-to-controller and controller-to-actuator communication paths, the control input is held constant between updates via a zero-order-hold (ZOH) mechanism: the value computed at the last triggering instant $t_\kappa$ is applied continuously over the interval $[t_\kappa,t_{\kappa+1})$, such that
\begin{align}
u(t)=K\xi(t_{\kappa})\,, \quad \forall t \!\in\! \lbrack t_{\kappa}\,, t_{\kappa+1}\phantom{(}\!\!) \,, \quad \kappa\!\in\! \mathbb{N}  \,. \label{eq:u_ET_v1}
\end{align}
The filters (\ref{eq:F_eta})--(\ref{eq:F_xi}) run continuously on the controller side, so that the filtered gradient estimate $\xi(t)$ is available at every instant, whereas only the correction $K\xi(t_{\kappa})$ is computed and transmitted at the event times. Since control updates occur at discrete time instants rather than continuously, we introduce the sampling-induced error
\begin{align}
e(t):=\xi(t_{\kappa})-\xi(t) \,, ~ \forall t \in \lbrack t_{\kappa}\,, t_{\kappa+1}\phantom{(}\!\!) \,,~ \kappa\in \mathbb{N} \,, \label{eq:e_ET_v1}
\end{align}
with which the event-triggered tuning law (\ref{eq:u_ET_v1}) is rewritten as
\begin{align}
u(t)&=K\xi(t)+Ke(t)  \,. \label{eq:u_ET_v2}
\end{align}
Plugging (\ref{eq:u_ET_v2}) into (\ref{eq:gradient_dtildeThetadt_20250206_1}) gives $\frac{d\tilde{p}}{dt}=(I_{2}+\Upsilon(t))K\xi+(I_{2}+\Upsilon(t))Ke$, while the filter equations (\ref{eq:F_dxi}) and (\ref{eq:F_deta}) do not depend on $u$ and remain unchanged. Therefore, with the state (\ref{eq:F_z}) and $B(t)$ defined in (\ref{eq:F_B}), the closed-loop system between two consecutive events is, exactly,
\begin{align}
\frac{dz}{dt}(t)=A(t)z(t)+b(t)+N(t,z(t))+B(t)e(t)\,, \label{eq:F_cl_et}
\end{align}
for all $t\in[t_{\kappa},t_{\kappa+1})$. The state $z(t)$ is continuous at every event: only the error $e(t)$ jumps, being reset to zero at $t=t_{\kappa}$. Moreover, (\ref{eq:F_cl_et}) is the continuous-time loop (\ref{eq:F_cl}) perturbed by $B(t)e(t)$, which is the situation covered by Lemma~\ref{lem:lyapunov}.

\subsection{Event-Triggering Mechanism}

Definition~\ref{def:staticEvent} shows how the error $e(t)$ --- the deviation between the last broadcasted value $\xi(t_{\kappa})$ of the filtered gradient estimate and its current value $\xi(t)$ --- is used to construct a static event-triggered mechanism for multivariable gradient-based source seeking. The mechanism combines a relative threshold, set by the design parameters $\sigma$ and $\alpha$, with an absolute threshold $\varepsilon_{0}$. It triggers the recalculation of the control law (\ref{eq:u_ET_v1}) and refreshes the ZOH actuator signal only when needed, as shown in the block diagram of Fig.~\ref{fig:BD_NewtonES_v5}, within the broader event-triggered sampling framework surveyed in \cite{HJT:2012}.

\begin{definition}[{\bf Event-Triggering Condition}] \label{def:staticEvent}
Let $\sigma\in(0,1)$, $\alpha>0$ and $\varepsilon_{0}>0$ be design parameters. The gradient-based source seeking controller with a static event-triggering condition consists of two components:
\begin{enumerate}
	\item A sequence of increasing times $\mathcal{I}=\{t_{0}\,, t_{1}\,, t_{2}\,,\ldots\}$, with $t_{0}=0$ and
				\begin{align}
					t_{\kappa+1}&=\inf\left\{t > t_{\kappa}:~ \alpha\|e(t)\|\geq \sigma\|\xi(t)\|+\varepsilon_{0} \right\}, \label{eq:tk+1_event}
				\end{align}
	where $\mathcal{I}=\{t_{0}\,, t_{1}\,, \ldots, t_{\kappa}\}$ if the set in (\ref{eq:tk+1_event}) is empty.
	\item A feedback control action, given by (\ref{eq:u_ET_v1}), updated at the time instants $t_{\kappa}$.
\end{enumerate}
\end{definition}

The two thresholds in (\ref{eq:tk+1_event}) play complementary roles. The relative term $\sigma\|\xi(t)\|$ tolerates large deviations far from the source, where the correction is large and need not be accurate. The absolute term $\varepsilon_{0}$ is needed because $\xi(t)$ becomes small as the vehicle approaches the source and may pass through zero; with a purely relative threshold, the triggering condition would then be met by arbitrarily small errors and no positive lower bound on the inter-event times could be guaranteed. With $\varepsilon_{0}>0$, the error $e(t)$ must grow from zero to at least $\varepsilon_{0}/\alpha$ between two consecutive events, which takes a time that does not depend on the probing frequency, since the rate of change of $\xi(t)$ is bounded independently of $\omega_{0}$. The absolute threshold adds a term of order $\varepsilon_{0}$ to the ultimate bound of Theorem~\ref{thm:et_ssc}.

The parameter $\sigma\in(0,1)$ is the fraction by which the decay rate of the continuous-time loop is reduced to allow less frequent updates, whereas the gain $\alpha$ and the threshold $\varepsilon_{0}$ must satisfy $\alpha\geq\alpha^{\ast}$ and $\varepsilon_{0}\leq\bar{\varepsilon}_{0}$, with $\alpha^{\ast}$ and $\bar{\varepsilon}_{0}$ given in (\ref{eq:F_rstar}). These bounds do not require the exact values of the plant parameters: $H^{\ast}$, $p^{\ast}$ and $f^{\ast}$ remain unknown, and only rough prior bounds on $f_{x}$ and $f_{y}$ are needed. A conservative choice (large $\alpha$ and small $\varepsilon_{0}$) only makes the events more frequent, without affecting stability.

\subsection{Convergence and Stability Guarantees}\label{sec:stability_analysis}

The following theorem establishes the stability of the scheme of Fig.~\ref{fig:BD_NewtonES_v5} and a minimum dwell time independent of the probing frequency. Let $\sigma\in(0,1)$ be the triggering parameter of Definition~\ref{def:staticEvent} and let the analysis constant $\sigma_{0}$ be chosen such that $\sigma_{0}\in(0,1-\sigma)$.
In addition to the constants $\alpha_{1}$, $\alpha_{2}$, $C_{N}$, $C_{B}$, $C_{p}$, $\rho_{0}$, $\rho$ and $\omega_{0}^{\ast}$ defined in (\ref{eq:F_alphas}), (\ref{eq:F_CN}), (\ref{eq:F_B}), (\ref{eq:F_Cp}), (\ref{eq:F_dbar}) and (\ref{eq:F_w0star}), define
\begin{align}
\alpha^{\ast}&:=2\alpha_{2}C_{B}\,, \quad \bar{\varepsilon}_{0}:=\frac{\rho_{0}}{2}\sqrt{\frac{\alpha_{1}\sigma_{0}(1-\sigma-\sigma_{0})}{2\alpha_{2}}}\,, \label{eq:F_rstar}\\
\bar{\omega}_{0}&:=\max\Big\{\omega_{0}^{\ast}\,,\frac{C_{p}}{\bar{\varepsilon}_{0}}\Big\}\,, \quad \rho_{D}:=\rho_{0}+\frac{C_{p}}{\bar{\omega}_{0}}+1\,, \label{eq:F_w0bar}\\
F_{D}&:=\frac{\|H^{\ast}\|}{2}\rho_{D}^{2}+a_{\mu}\|H^{\ast}\|\rho_{D}+\|H^{\ast}\|a_{\mu}^{2}+\rho_{D}\,, \label{eq:F_FD}\\
L_{D}&:=\omega_{l}\big(\rho_{D}+\mu\, m\, F_{D}\big)\,. \label{eq:F_LD}
\end{align}
All these constants are independent of $\omega_{0}$.

\begin{theorem} \label{thm:et_ssc}
Consider the closed-loop system (\ref{eq:F_cl_et}), composed of the vehicle (\ref{eq:dotX_ETSSC_v1})--(\ref{eq:dotTheta_ETSSC_v1}) with the velocities (\ref{eq:vt_ETSSC_v1})--(\ref{eq:wt_ETSSC_v1}), the filters (\ref{eq:F_eta})--(\ref{eq:F_xi}), the event-triggered tuning law (\ref{eq:u_ET_v1}) and the event-triggering mechanism of Definition~\ref{def:staticEvent}, under Assumptions~\ref{assum:f} and~\ref{assum:w}. For all $\alpha\geq\alpha^{\ast}$, $\varepsilon_{0}\in(0,\bar{\varepsilon}_{0}]$, $\omega_{0}\geq\bar{\omega}_{0}$ and $\|z(0)\|\leq\rho$, the following holds:
\begin{enumerate}
\item[(a)] The closed-loop solution exists for all $t\geq 0$ and is unique.
\item[(b)] The inter-event times satisfy
\begin{align}
t_{\kappa+1}-t_{\kappa}\geq\tau^{\ast}:=\frac{\varepsilon_{0}}{\alpha L_{D}}>0\,, \quad \forall\kappa\in\mathbb{N}\,, \label{eq:F_dwell}
\end{align}
where $L_{D}$ is independent of $\omega_{0}$; hence Zeno behavior is excluded.
\item[(c)] The vehicle position satisfies, for all $t\geq 0$,
\begin{align}
\|p(t)-p^{\ast}\|&\leq c_{1}e^{-\lambda(\sigma) t}\|z(0)-z_{\rm per}(0)\| \nonumber \\
&\quad+\mathcal{O}\left(\frac{\mu^{2}}{\mu^{2}-1}a_{1}+\varepsilon_{0}+\frac{1}{\omega_{0}}\right), \label{eq:norm_p_thm2}
\end{align}
where $c_{1}:=\sqrt{\alpha_{2}/\alpha_{1}}$ and the decay rate depends on the triggering parameter $\sigma$ as
\begin{align}
\lambda(\sigma):=\frac{1}{2}\frac{1}{\alpha_{2}}\big(1-\sigma-\sigma_{0}\big)>0\,, \label{eq:F_lambda}
\end{align}
where $z_{\rm per}(t)$ is the periodic solution of the continuous-time loop given by Theorem~\ref{thm:ssc} and the constant implicit in $\mathcal{O}(\cdot)$ is independent of $\omega_{0}$ and $\varepsilon_{0}$.
\end{enumerate}
\end{theorem}

\begin{proof}
Between two consecutive events the held value $\xi(t_{\kappa})$ is constant and the right-hand side of (\ref{eq:F_cl_et}), with $e=\xi(t_{\kappa})-\xi$, is continuous in $t$ and polynomial in $z$. Hence, from each event time there is a unique solution up to the next event time, and the concatenation of these pieces defines a unique solution on a maximal interval $[0,t_{\max})$, where $t_{\max}\leq\infty$ is either a finite escape time or the accumulation point of the event times. Since $\omega_{0}\geq\bar{\omega}_{0}\geq\omega_{0}^{\ast}$, Lemmas~\ref{lem:periodic} and~\ref{lem:lyapunov} apply; let $\tilde{z}(t):=z(t)-z_{\rm per}(t)$ and $V(t):=\tilde{z}^{\top}(t)P(t)\tilde{z}(t)$. Because $z(t)$ is continuous at the events, $\tilde{z}(t)$ and $V(t)$ are continuous on $[0,t_{\max})$, and they are continuously differentiable between events. The proof is divided into two parts: stability analysis and avoidance of Zeno behavior.

\emph{Part A (stability analysis). Step A1: Lyapunov inequality.} By Lemma~\ref{lem:lyapunov}, inequality (\ref{eq:F_Vdot}) holds on every inter-event interval $(t_{\kappa},t_{\kappa+1})$.

\emph{Step A2: bound on the sampling-induced error.} Since $e(t_{\kappa})=0$ and $t_{\kappa+1}$ is the first instant after $t_{\kappa}$ at which the condition in (\ref{eq:tk+1_event}) is met, $\alpha\|e(t)\|<\sigma\|\xi(t)\|+\varepsilon_{0}$ for all $t\in[t_{\kappa},t_{\kappa+1})$. Moreover, $\|\xi\|\leq\|z\|\leq\|\tilde{z}\|+\|z_{\rm per}\|\leq\|\tilde{z}\|+C_{p}/\omega_{0}$ by (\ref{eq:F_zper_bound}). Hence, using $\sigma<1$,
\begin{align}
\alpha\|e(t)\|\leq \sigma\|\tilde{z}(t)\|+\delta\,, \quad \delta:=\varepsilon_{0}+\frac{C_{p}}{\omega_{0}}\,, \label{eq:F_ebound}
\end{align}
for all $t\in[t_{\kappa},t_{\kappa+1})$ and all $\kappa$.

\emph{Step A3: decay on the ball $\|\tilde{z}\|\leq\rho_{0}$.} Since $\alpha\geq\alpha^{\ast}=2\alpha_{2}C_{B}$ by (\ref{eq:F_rstar}), the last term of (\ref{eq:F_Vdot}) satisfies $2\alpha_{2}C_{B}\|\tilde{z}\|\|e\|=\frac{2\alpha_{2}C_{B}}{\alpha}\|\tilde{z}\|\,\alpha\|e\|\leq\|\tilde{z}\|\,\alpha\|e\|$. Substituting (\ref{eq:F_ebound}),
\begin{align}
\frac{dV}{dt}&\leq-\|\tilde{z}\|^{2}+2\alpha_{2}C_{N}\|\tilde{z}\|^{3}+\sigma\|\tilde{z}\|^{2}+\delta\|\tilde{z}\|\,. \label{eq:F_VA3}
\end{align}
The three positive terms of (\ref{eq:F_VA3}) are treated as follows. (i) The term $\sigma\|\tilde{z}\|^{2}$, generated by the relative threshold, is quadratic in $\|\tilde{z}\|$, like the negative term, and this is what produces the factor $1-\sigma$ in the decay rate. (ii) If $\|\tilde{z}\|\leq\rho_{0}$, the cubic term satisfies $2\alpha_{2}C_{N}\|\tilde{z}\|^{3}\leq\frac{\sigma_{0}}{2}\|\tilde{z}\|^{2}$ by the definition of $\rho_{0}$ in (\ref{eq:F_dbar}). (iii) The last term, generated by the absolute threshold $\varepsilon_{0}$ and by the periodic solution $z_{\rm per}$, is linear in $\|\tilde{z}\|$ and cannot be dominated by $-\|\tilde{z}\|^{2}$ near $\tilde{z}=0$; by Young's inequality $2ab\leq\frac{\sigma_{0}}{2}a^{2}+\frac{2}{\sigma_{0}}b^{2}$ with $a=\|\tilde{z}\|$ and $b=\frac{\delta}{2}$, $\delta\|\tilde{z}\|\leq\frac{\sigma_{0}}{2}\|\tilde{z}\|^{2}+\frac{\delta^{2}}{2\sigma_{0}}$. Therefore, on the ball $\|\tilde{z}\|\leq\rho_{0}$, using (\ref{eq:F_Pbounds}),
\begin{align}
\frac{dV}{dt}&\leq-\Big(1-\sigma-\frac{\sigma_{0}}{2}-\frac{\sigma_{0}}{2}\Big)\|\tilde{z}\|^{2}+\frac{\delta^{2}}{2\sigma_{0}} \nonumber \\
&=-\big(1-\sigma-\sigma_{0}\big)\|\tilde{z}\|^{2}+\frac{\delta^{2}}{2\sigma_{0}} \nonumber \\
&\leq-2\lambda(\sigma)V+\frac{\delta^{2}}{2\sigma_{0}}\,, \label{eq:F_VA3b}
\end{align}
with $\lambda(\sigma)$ defined in (\ref{eq:F_lambda}); the last step is valid because $1-\sigma-\sigma_{0}>0$, so that $-(1-\sigma-\sigma_{0})\|\tilde{z}\|^{2}\leq-\frac{1-\sigma-\sigma_{0}}{\alpha_{2}}V=-2\lambda(\sigma)V$.

\emph{Step A4: integration across the events.} Let $t^{\prime}\in(0,t_{\max}]$ be such that $\|\tilde{z}(t)\|\leq\rho_{0}$ for all $t\in[0,t^{\prime})$ and every compact subinterval of $[0,t^{\prime})$ contains finitely many events. Define
\begin{align}
V_{\infty}&:=\frac{\delta^{2}}{2\sigma_{0}\,2\lambda(\sigma)}=\frac{\alpha_{2}\delta^{2}}{2\sigma_{0}(1-\sigma-\sigma_{0})}\,, \label{eq:F_Vinfdef}\\
Y(t)&:=e^{2\lambda(\sigma)t}\big(V(t)-V_{\infty}\big)\,. \nonumber
\end{align}
The function $Y$ is continuous on $[0,t^{\prime})$ and, by (\ref{eq:F_VA3b}), $\frac{dY}{dt}=e^{2\lambda(\sigma)t}\big(\frac{dV}{dt}+2\lambda(\sigma)V-\frac{\delta^{2}}{2\sigma_{0}}\big)\leq 0$ between events. Hence $Y$ is nonincreasing on each inter-event interval and, being continuous at the events, on the whole interval $[0,t^{\prime})$; no product of estimates over the successive inter-event intervals is needed. From $Y(t)\leq Y(0)$,
\begin{align}
V(t)\leq e^{-2\lambda(\sigma)t}V(0)+\Big(1-e^{-2\lambda(\sigma)t}\Big)V_{\infty}\,, \label{eq:F_VA4}
\end{align}
for all $t\in[0,t^{\prime})$.

\emph{Part B (avoidance of Zeno behavior). Step B1: rate of change of the filtered gradient estimate on a compact set.} Let $D:=\{z\in\mathbb{R}^{5}:\ \|z\|\leq\rho_{D}\}$, with $\rho_{D}$ in (\ref{eq:F_w0bar}). For $z\in D$, (\ref{eq:F_hp}), Lemma~\ref{lem:dither} and $|q(\tilde{p})|\leq\frac{1}{2}\|H^{\ast}\|\|\tilde{p}\|^{2}$ give $|f-\eta|\leq|q(\tilde{p})|+\|S\|\|H^{\ast}\|\|\tilde{p}\|+|s-s_{0}|+|\tilde{\eta}| \leq\frac{\|H^{\ast}\|}{2}\rho_{D}^{2}+a_{\mu}\|H^{\ast}\|\rho_{D}+\|H^{\ast}\|a_{\mu}^{2}+\rho_{D}=F_{D}$, so that, by (\ref{eq:hatG_ETSSC_v1}) and (\ref{eq:F_Mbound}), the gradient estimate satisfies $\|\hat{G}\|=\|M\||f-\eta|\leq\mu\, m\, F_{D}$. Then (\ref{eq:F_xi}) gives
\begin{align}
\Big\|\frac{d\xi}{dt}\Big\|&\leq\omega_{l}\big(\|\xi\|+\|\hat{G}\|\big) \nonumber \\
&\leq\omega_{l}\big(\rho_{D}+\mu\, m\, F_{D}\big)=L_{D}\,. \label{eq:F_dxi_bound}
\end{align}
The constant $L_{D}$ is independent of $\omega_{0}$ because the right-hand side of (\ref{eq:F_xi}) contains $\hat{G}(t)$, and hence $M(t)$, whose amplitude is independent of $\omega_{0}$, but not the derivative of $M(t)$, whose amplitude grows linearly with $\omega_{0}$.

\emph{Step B2: dwell time on $D$.} Suppose that $z(t)\in D$ for all $t\in[t_{\kappa},t_{\kappa+1}]$. On this interval the held value $\xi(t_{\kappa})$ is constant, so $e(t_{\kappa})=0$ and $\frac{de}{dt}=-\frac{d\xi}{dt}$; thus, by (\ref{eq:F_dxi_bound}),
\begin{align}
\|e(t)\|=\Big\|\int_{t_{\kappa}}^{t}\frac{d\xi}{dt}(\tau)d\tau\Big\|\leq L_{D}(t-t_{\kappa})\,. \label{eq:F_egrowth}
\end{align}
On the other hand, by the definition of $t_{\kappa+1}$ as an infimum and the continuity of $\xi$, the triggering condition holds in the limit $t\to t_{\kappa+1}$ from the left: $\alpha\|\xi(t_{\kappa})-\xi(t_{\kappa+1})\|\geq \sigma\|\xi(t_{\kappa+1})\|+\varepsilon_{0}\geq\varepsilon_{0}$. Combining with (\ref{eq:F_egrowth}), $\varepsilon_{0}\leq\alpha L_{D}(t_{\kappa+1}-t_{\kappa})$, i.e.,
\begin{align}
t_{\kappa+1}-t_{\kappa}\geq\frac{\varepsilon_{0}}{\alpha L_{D}}\,. \label{eq:F_dwellB2}
\end{align}
This step uses no stability estimate, so there is no circularity with Part~A.

\emph{Step B3: invariance and continuation.} Since $\|z(0)\|\leq\rho$ and $\|z_{\rm per}(0)\|\leq C_{p}/\omega_{0}\leq\rho$, one has $\|\tilde{z}(0)\|\leq 2\rho=\frac{\rho_{0}}{2}\sqrt{\alpha_{1}/\alpha_{2}}\leq\frac{\rho_{0}}{2}$ and, by (\ref{eq:F_Pbounds}),
\begin{align}
V(0)\leq\alpha_{2}\|\tilde{z}(0)\|^{2}\leq 4\alpha_{2}\rho^{2}=\frac{\alpha_{1}\rho_{0}^{2}}{4}\,. \label{eq:F_V0}
\end{align}
Moreover, $\varepsilon_{0}\leq\bar{\varepsilon}_{0}$ and $C_{p}/\omega_{0}\leq C_{p}/\bar{\omega}_{0}\leq\bar{\varepsilon}_{0}$ give $\delta\leq 2\bar{\varepsilon}_{0}$, so that, by the definition of $\bar{\varepsilon}_{0}$ in (\ref{eq:F_rstar}),
\begin{align}
V_{\infty}\leq\frac{2\alpha_{2}\bar{\varepsilon}_{0}^{2}}{\sigma_{0}(1-\sigma-\sigma_{0})}=\frac{\alpha_{1}\rho_{0}^{2}}{4}\,. \label{eq:F_Vinf}
\end{align}
Let $t_{1}:=\sup\{t\in[0,t_{\max}):\ \|\tilde{z}(\tau)\|<\rho_{0}~~\forall\tau\in[0,t]\}>0$. For $t\in[0,t_{1})$, $\|z(t)\|\leq\|\tilde{z}(t)\|+\|z_{\rm per}(t)\|<\rho_{0}+C_{p}/\bar{\omega}_{0}<\rho_{D}$, so $z(t)\in D$; by Step~B2, any two consecutive events in $[0,t_{1})$ are separated by at least $\varepsilon_{0}/(\alpha L_{D})$, and every compact subinterval of $[0,t_{1})$ contains finitely many events. Thus Step~A4 applies with $t^{\prime}=t_{1}$, and the right-hand side of (\ref{eq:F_VA4}), being a convex combination of $V(0)$ and $V_{\infty}$, is bounded by $\max\{V(0),V_{\infty}\}$. From (\ref{eq:F_V0}), (\ref{eq:F_Vinf}) and (\ref{eq:F_Pbounds}),
\begin{align}
\alpha_{1}\|\tilde{z}(t)\|^{2}\leq V(t)\leq\frac{\alpha_{1}\rho_{0}^{2}}{4}\,, \quad \forall t\in[0,t_{1})\,, \label{eq:F_inv}
\end{align}
i.e., $\|\tilde{z}(t)\|\leq\rho_{0}/2$ on $[0,t_{1})$.
If $t_{1}<t_{\max}$, the continuity of $\tilde{z}$ would give $\|\tilde{z}(t_{1})\|\leq\rho_{0}/2<\rho_{0}$, and $\|\tilde{z}\|$ would remain below $\rho_{0}$ on a right neighborhood of $t_{1}$, contradicting the definition of $t_{1}$; hence $t_{1}=t_{\max}$. If $t_{\max}<\infty$, then on $[0,t_{\max})$ the state $z(t)$ remains in the compact set $D$, which excludes a finite escape time, and the events are separated by at least $\varepsilon_{0}/(\alpha L_{D})$, which excludes an accumulation of event times; the solution could then be continued beyond $t_{\max}$, a contradiction. Therefore $t_{\max}=\infty$, which proves (a); $z(t)\in D$ for all $t\geq 0$, so (\ref{eq:F_dwellB2}) holds for every $\kappa$, which proves (b); and the sublevel set $\{V\leq\alpha_{1}\rho_{0}^{2}/4\}$ is forward invariant.

\emph{Conclusion of the proof of (c).} By Step~B3, (\ref{eq:F_VA4}) holds for all $t\geq 0$. Using (\ref{eq:F_Pbounds}) and $\sqrt{a+b}\leq\sqrt{a}+\sqrt{b}$,
\begin{align}
\|\tilde{z}(t)\|&\leq\sqrt{\frac{\alpha_{2}}{\alpha_{1}}}\,e^{-\lambda(\sigma)t}\|\tilde{z}(0)\|+\sqrt{\frac{V_{\infty}}{\alpha_{1}}} \nonumber \\
&= c_{1}e^{-\lambda(\sigma) t}\|\tilde{z}(0)\|+c_{3}\Big(\varepsilon_{0}+\frac{C_{p}}{\omega_{0}}\Big)\,, \label{eq:F_d_et}\\
c_{3}&:=\frac{c_{1}}{\sqrt{2\sigma_{0}(1-\sigma-\sigma_{0})}}\,, \nonumber
\end{align}
where the equality follows from (\ref{eq:F_Vinfdef}) and $c_{1}=\sqrt{\alpha_{2}/\alpha_{1}}$. Finally, from (\ref{eq:q_ETSSC_v2}), $\|\tilde{p}\|\leq\|z\|\leq\|\tilde{z}\|+\|z_{\rm per}\|$, (\ref{eq:F_Sbound}), (\ref{eq:F_zper_bound}) and (\ref{eq:F_d_et}), one can obtain

\begin{align}
\|p(t)-p^{\ast}\|&\leq\|\tilde{p}(t)\|+\|S(t)\|\leq\|\tilde{z}(t)\|+\frac{C_{p}}{\omega_{0}}+a_{\mu} \nonumber \\
&\leq c_{1}e^{-\lambda(\sigma) t}\|z(0)-z_{\rm per}(0)\|+\frac{\mu^{2}}{\mu^{2}-1}a_{1} \nonumber \\
&\quad+c_{3}\varepsilon_{0}+\frac{(c_{3}+1)C_{p}}{\omega_{0}}\,. \label{eq:F_p_et}
\end{align}
Since the coefficients $c_{3}$ and $(c_{3}+1)C_{p}$ are independent of $\omega_{0}$ and $\varepsilon_{0}$, the last three terms of (\ref{eq:F_p_et}) are bounded by $c_{2}\big(\frac{\mu^{2}}{\mu^{2}-1}a_{1}+\varepsilon_{0}+\frac{1}{\omega_{0}}\big)$, with $c_{2}:=\max\{1\,,c_{3}\,,(c_{3}+1)C_{p}\}$, which is (\ref{eq:norm_p_thm2}).
\end{proof}

\begin{remark}[Influence of the triggering parameters] \label{rem:trigger_params}
Since $P(t)$ satisfies (\ref{eq:F_Pdot}) with $Q=I_{5}$ and $\alpha_{2}\geq\lambda_{\max}(P(t))$, the decay rate (\ref{eq:F_lambda}) has the classical form $\frac{1}{2}\frac{\lambda_{\min}(Q)}{\lambda_{\max}(P)}(1-\sigma-\sigma_{0})$. It decreases linearly with the relative threshold $\sigma$, from the rate of the continuous-time law in Theorem~\ref{thm:ssc}, recovered as $\sigma\to 0$, down to zero as $\sigma\to 1-\sigma_{0}$: a larger tolerance on the sampling-induced error, which yields fewer updates, comes with slower convergence, as in static relative triggering mechanisms \cite{AUT:2025}. The constant $\sigma_{0}$ cannot be removed, since it accounts for the two terms of (\ref{eq:F_VA3}) that are not quadratic in the error, but it can be chosen arbitrarily small in $(0,1-\sigma)$, which shrinks the region of attraction, since $\rho$ in (\ref{eq:F_dbar}) is proportional to $\sigma_{0}$, and enlarges the constant $c_{3}$ in the ultimate bound. The gain $\alpha$ does not appear in the rate or in the ultimate bound; a larger $\alpha$ only makes the events more frequent. The absolute threshold $\varepsilon_{0}$ does not affect the rate either: it enters the ultimate bound, through the term of order $\varepsilon_{0}$ in (\ref{eq:norm_p_thm2}), and the dwell time $\tau^{\ast}=\varepsilon_{0}/(\alpha L_{D})$, which is independent of $\omega_{0}$ and therefore does not shrink when the probing frequency is increased to tighten the term of order $1/\omega_{0}$. Finally, as in Remark~\ref{rem:init}, (\ref{eq:norm_p_thm2}) implies the convergence of the measured signal to a neighborhood of $f^{\ast}$ of order $\mathcal{O}\big(\frac{\mu^{4}}{(\mu^{2}-1)^{2}}a_{1}^{2}+\varepsilon_{0}^{2}+\frac{1}{\omega_{0}^{2}}\big)$.
\end{remark}

\section{SIMULATION RESULTS}
In the following simulations, we choose $\mu = 5$ and $\omega_{0} = 20$ rad/s, so that, by Assumption~\ref{assum:w}, $\omega_{1} = \mu\,\omega_{0} = 100$ rad/s, and the probing period is $T=2\pi/\omega_{0}\approx 0.314$ s. As in \cite{SCL_2007}, the dither amplitude is $a_{1} = 0.1$. The curvatures of the unknown map~\eqref{eq:Q_ETSSC_v1} are $f_x = 1$ and $f_y = 0.5$, with maximum value $f^{*} = 1$. The gain matrix is $K=\mathrm{diag}(1,1)$. The cut-off frequencies of the high-pass and low-pass filters \eqref{eq:F_eta} and \eqref{eq:F_xi} are $\omega_{h}=0.75$ rad/s and $\omega_{l}=1.3$ rad/s, and both filters are initialized at rest, $\eta(0)=0$ and $\xi(0)=0$ (see Remark~\ref{rem:init}). For these gains, the eigenvalues of the average system matrix $\bar{A}$ in \eqref{eq:F_Abar} are $-\omega_{h}=-0.75$ and the roots of \eqref{eq:F_charpoly}, namely the roots $-0.65\pm 0.937j$ of $\lambda^{2}+1.3\lambda+1.3=0$ and the roots $-0.65\pm 0.477j$ of $\lambda^{2}+1.3\lambda+0.65=0$. The event-triggering parameters in Definition~\ref{def:staticEvent} are $\sigma=0.75$, $\alpha=1$ and $\varepsilon_{0}=5\times 10^{-2}$. The closed-loop equations are integrated with the forward Euler method with a fixed step of $10^{-3}$ s. The controller analyzed in Sections~\ref{sec:CT} and~\ref{sec:ETSSC} and the controller simulated here are identical. The thresholds $\omega_{0}^{\ast}$, $\bar{\omega}_{0}$, $\alpha^{\ast}$ and $\bar{\varepsilon}_{0}$ of Theorems~\ref{thm:ssc} and~\ref{thm:et_ssc} are sufficient conditions derived from worst-case estimates and are conservative; the simulations use considerably less restrictive values of $\omega_{0}$, $\alpha$ and $\varepsilon_{0}$. Moreover, Theorems~\ref{thm:ssc} and~\ref{thm:et_ssc} assume a fixed source: for a mobile source, the dynamics of $\tilde{p}(t)$ acquire the additional term $-\dot{p}^{\ast}(t)$, which is not covered by the analysis, so that the mobile-source simulations are included only as a numerical illustration of the behavior of the same controller for a slowly moving source.

\begin{figure}[h!]
	\centering
	\subfigure[Vehicle Trajectory, $(x,y)$. \label{fig:SSC_S_xy}]
	{\includegraphics[width=4.26cm]{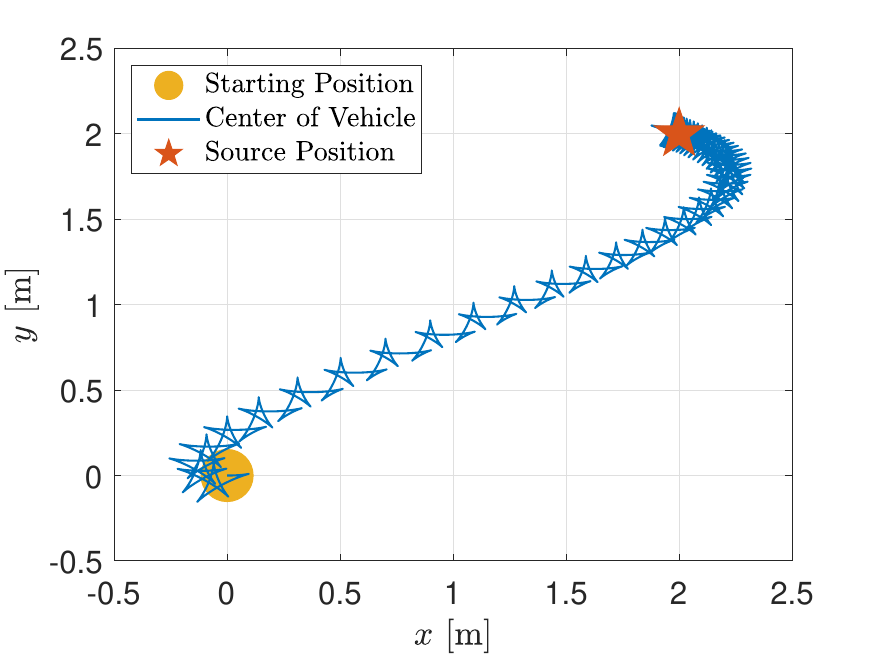}}
	\subfigure[Tuning Law, $u(t)$. \label{fig:SSC_S_u}]
	{\includegraphics[width=4.26cm]{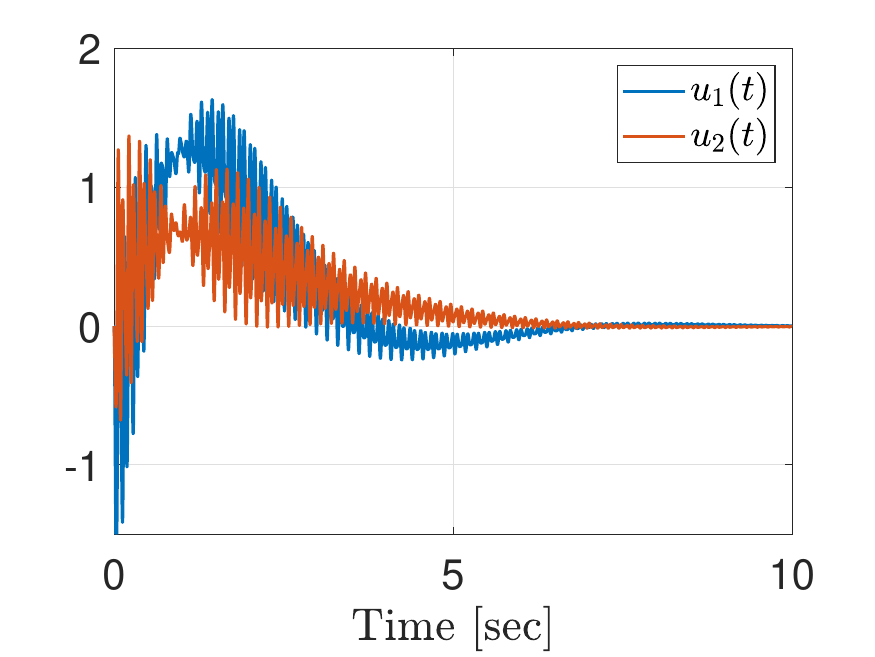}}
	\\
	\subfigure[Vehicle Trajectory, $(x,y)$. \label{fig:SSC_M_xy}]
	{\includegraphics[width=4.26cm]{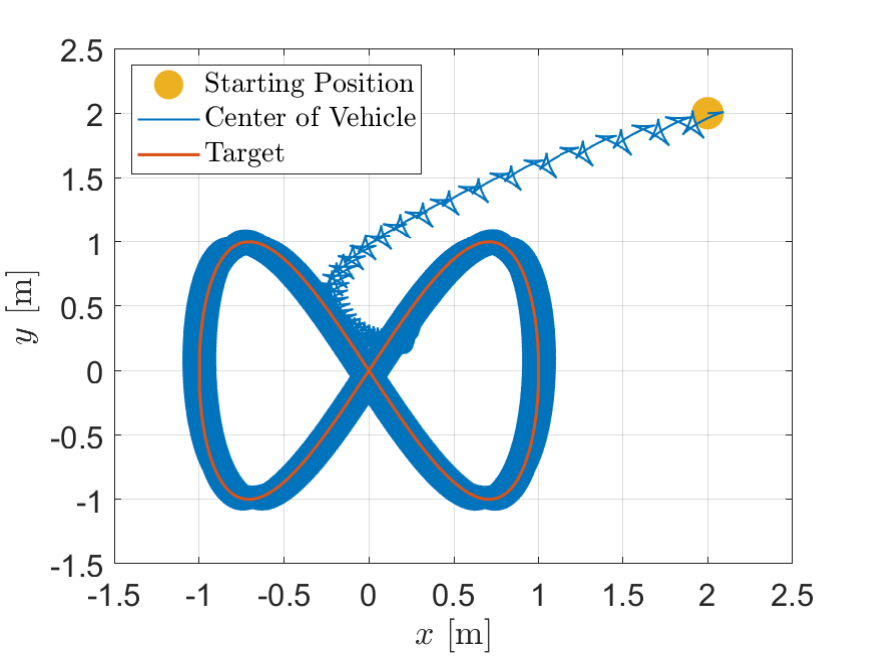}}
	\subfigure[Tuning Law, $u(t)$. \label{fig:SSC_M_u}]
	{\includegraphics[width=4.26cm]{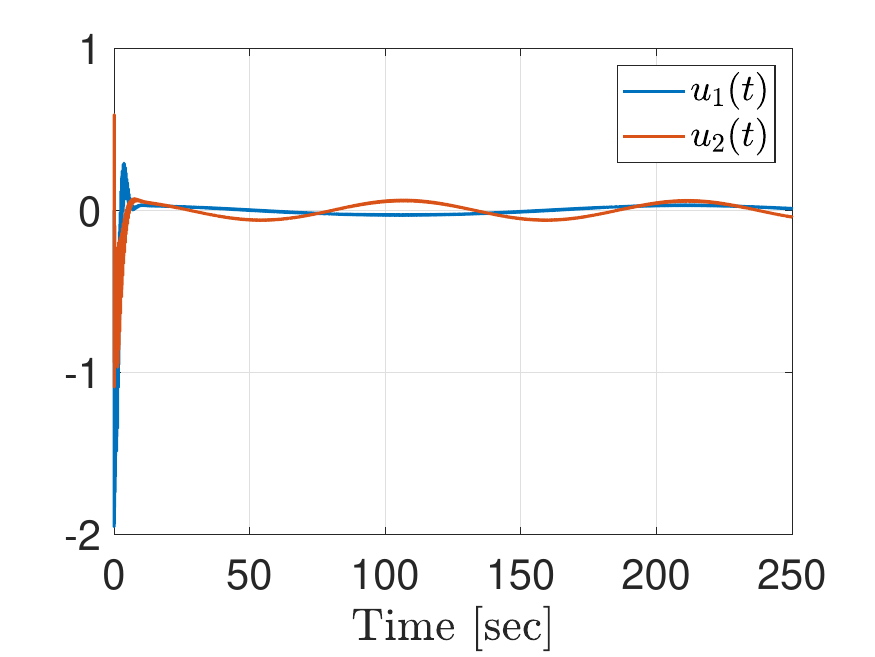}}
	\begin{picture}(0,0)
		\put(39,38){\includegraphics[width=1.85cm]{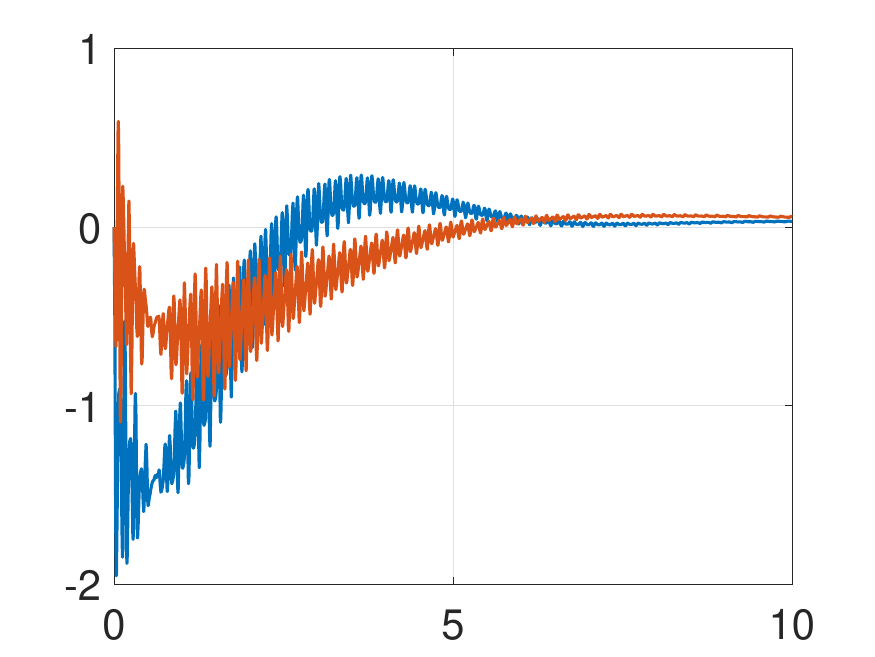}}
	\end{picture}
	\caption{Continuous Tuning Law: fixed and mobile source.}
	\label{fig:SSC}
\end{figure}

\begin{figure}[h!]
	\centering
	\subfigure[Output, $f(x,y)$. \label{fig:ET_SSC_S_f}]
	{\includegraphics[width=4.26cm]{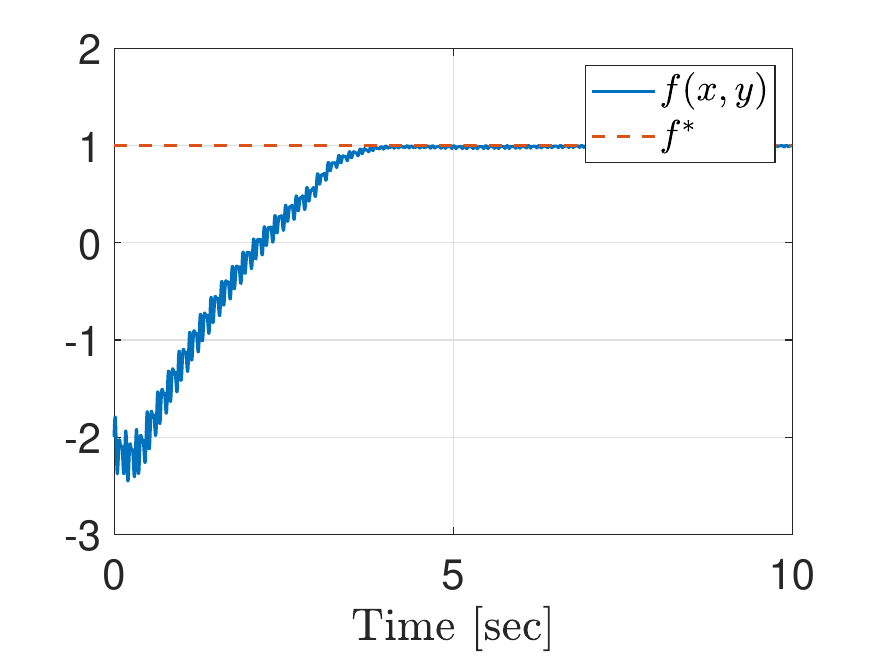}}
	\subfigure[Vehicle Position, $(x(t),y(t))$. \label{fig:ET_SSC_S_p}]
	{\includegraphics[width=4.26cm]{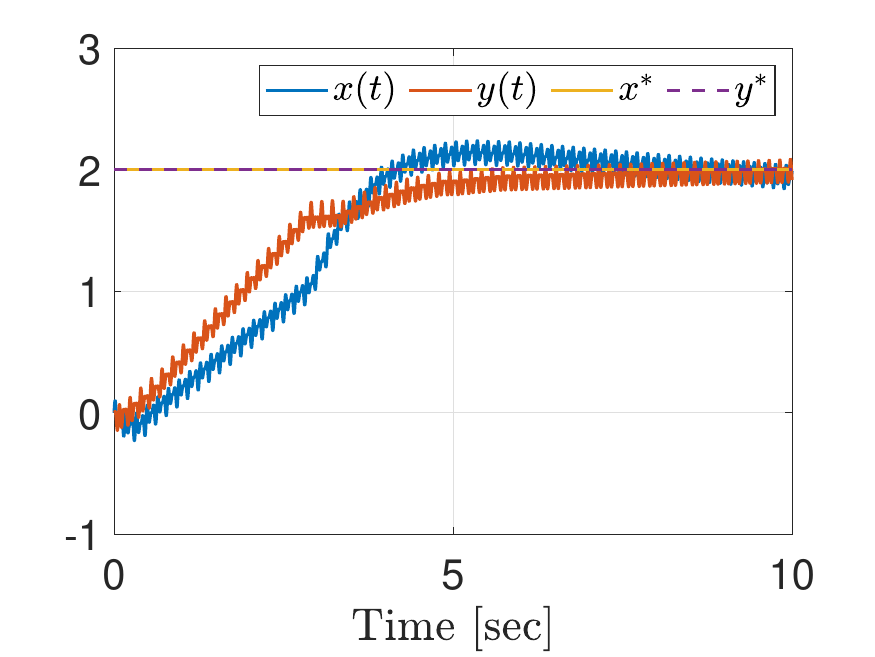}}
	\\
	\subfigure[Vehicle Trajectory, $(x,y)$. \label{fig:ET_SSC_S_xy}]
	{\includegraphics[width=4.26cm]{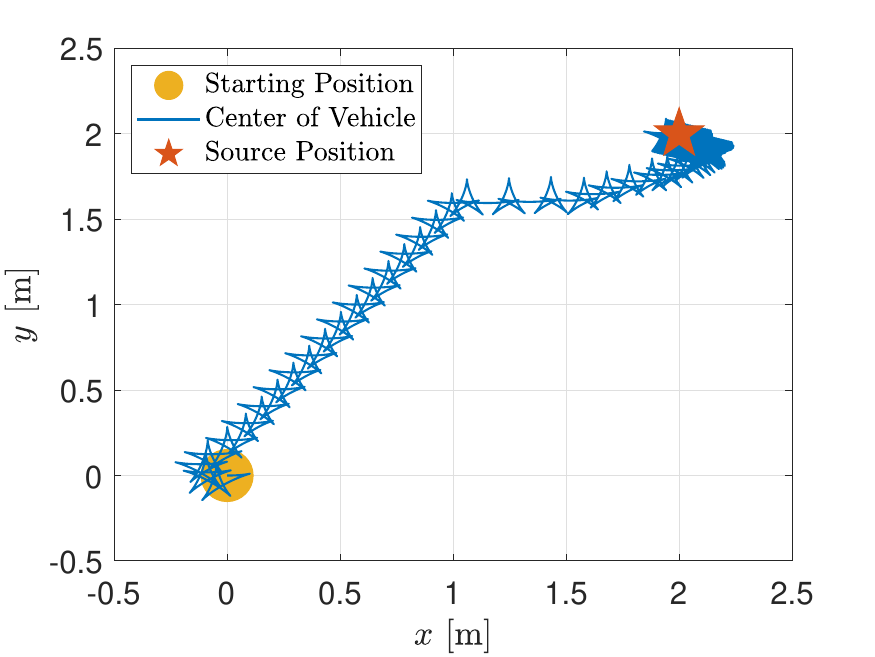}}
	\subfigure[Tuning Law, $u(t)$. \label{fig:ET_SSC_S_u}]
	{\includegraphics[width=4.26cm]{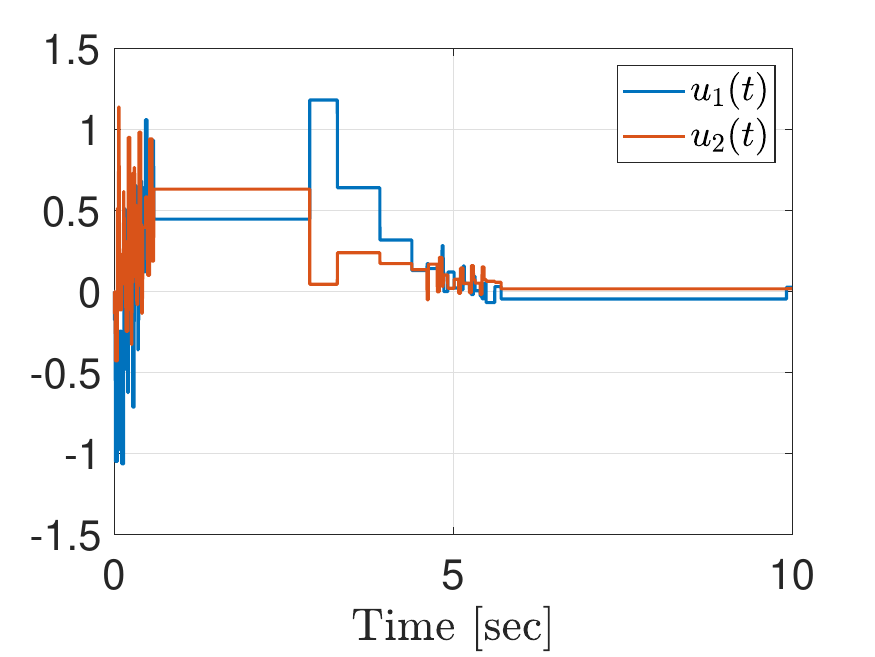}}
	\begin{picture}(0,0)
		\put(25,38){\includegraphics[width=1.5cm]{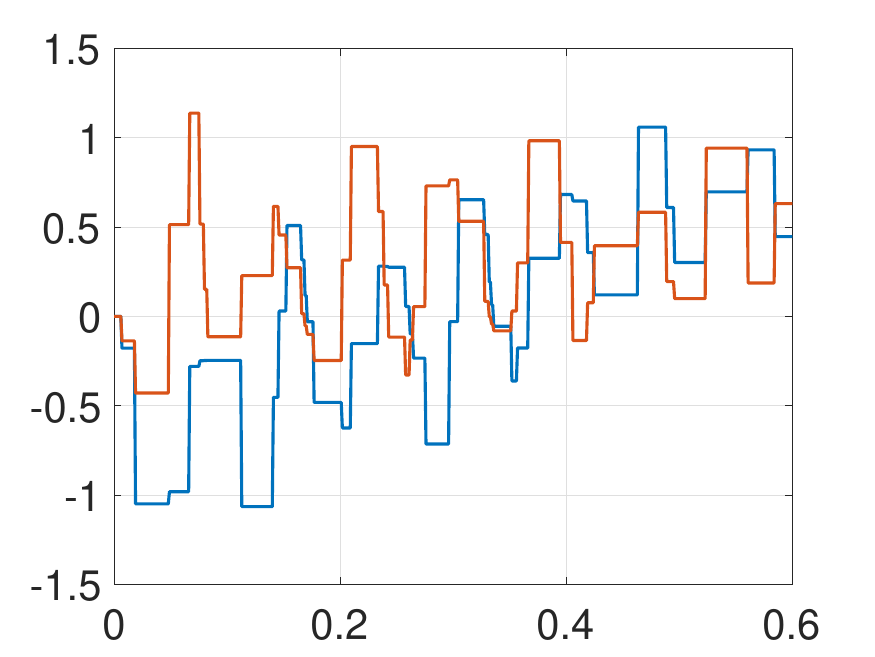}}
		\put(54,160){\includegraphics[width=2cm]{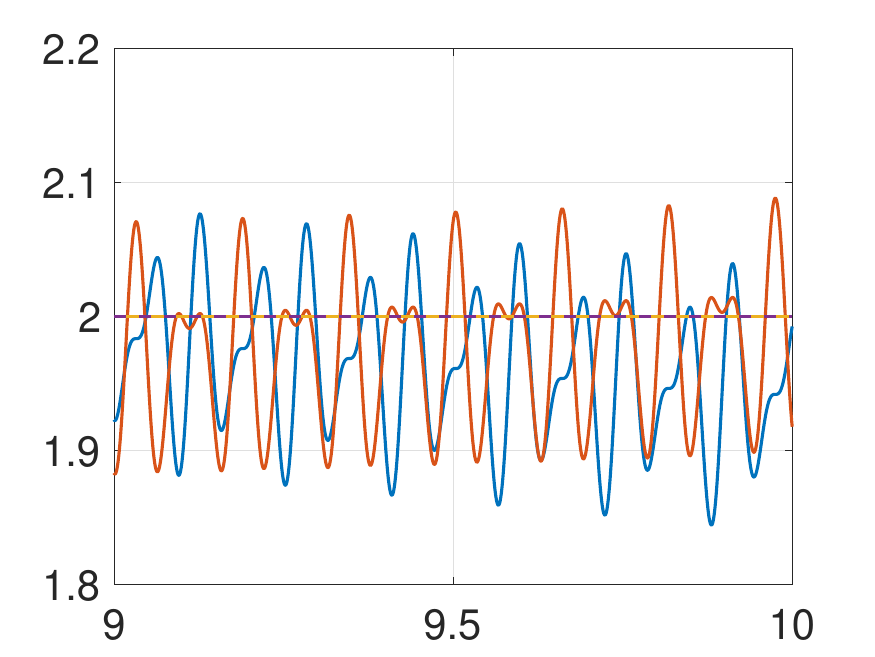}}
		\put(-78,160){\includegraphics[width=2.20cm]{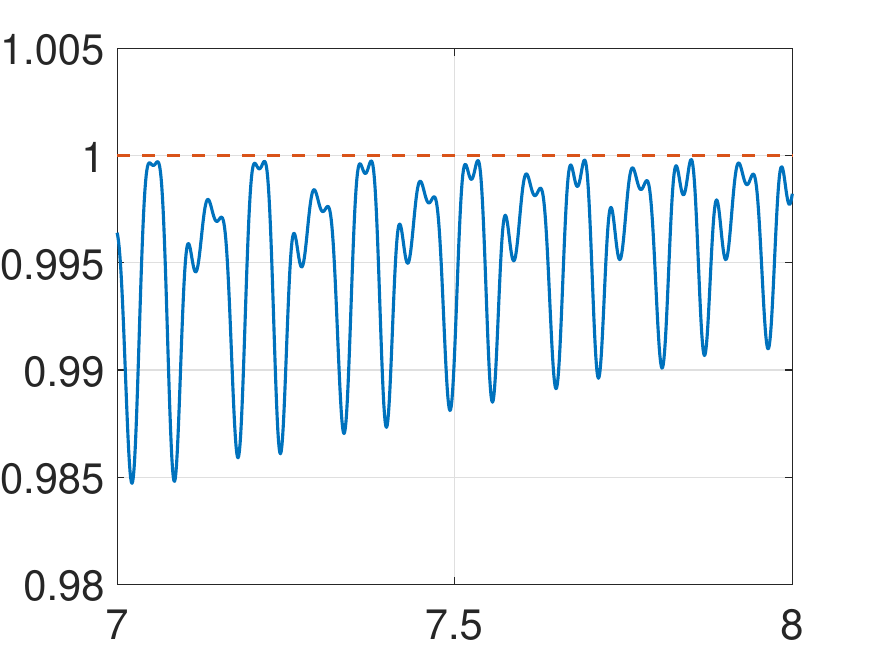}}
	\end{picture}
	\caption{Event-Triggered Tuning Law: fixed source.}
	\label{fig:ET_SSC_S}
\end{figure}

\begin{figure}[h!]
	\centering
	\subfigure[Output, $f(x,y)$. \label{fig:ET_SSC_M_f}]
	{\includegraphics[width=4.26cm]{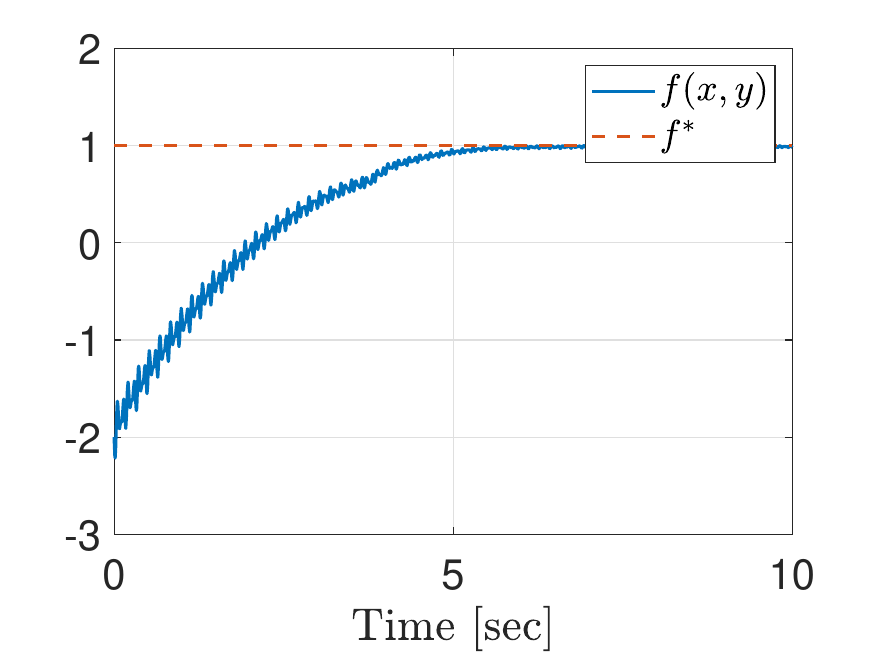}}
	\subfigure[Vehicle Position, $(x(t),y(t))$. \label{fig:ET_SSC_M_p}]
	{\includegraphics[width=4.26cm]{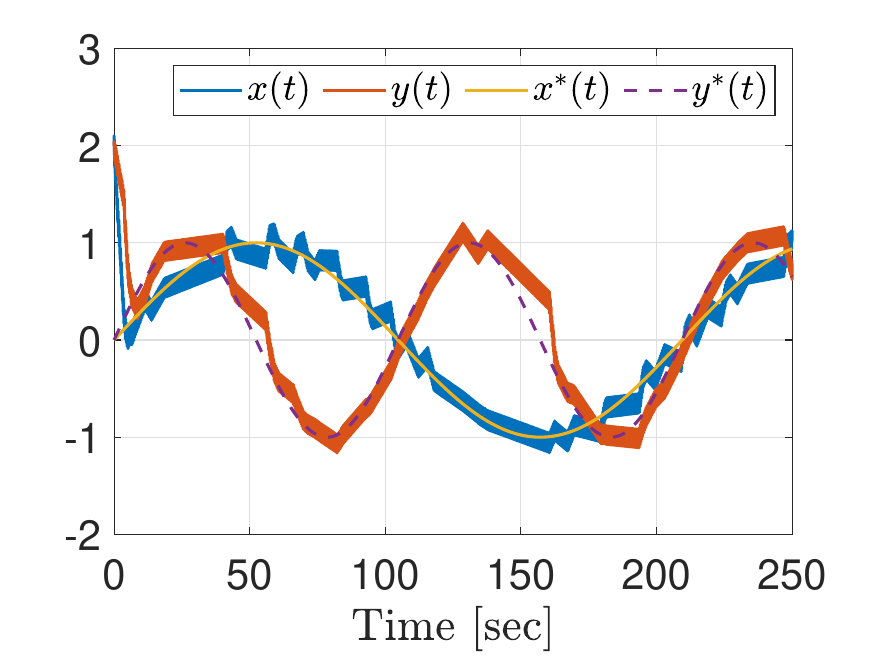}}
	\\
	\subfigure[Vehicle Trajectory, $(x,y)$. \label{fig:ET_SSC_M_xy}]
	{\includegraphics[width=4.26cm]{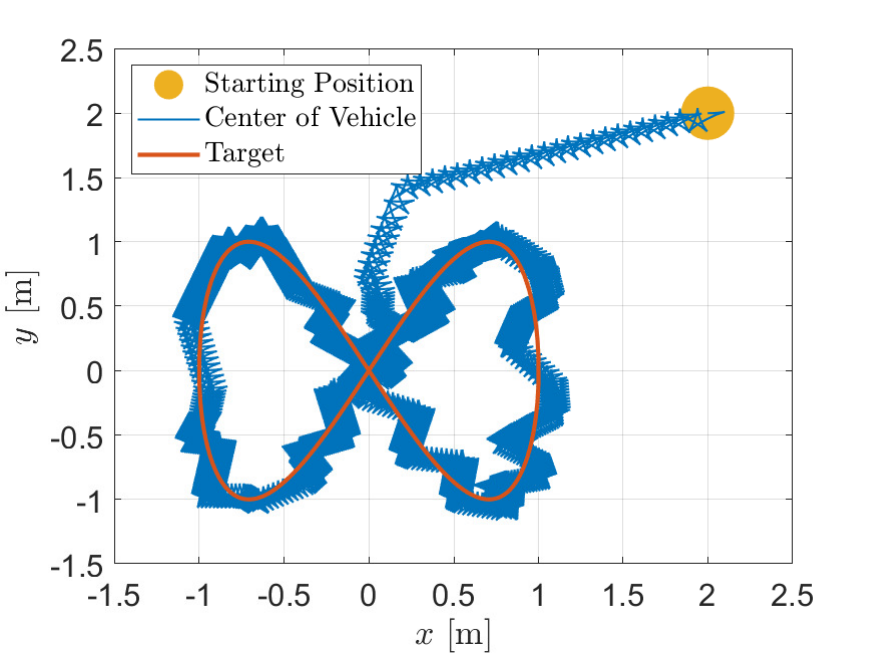}}
	\subfigure[Tuning Law, $u(t)$. \label{fig:ET_SSC_M_u}]
	{\includegraphics[width=4.26cm]{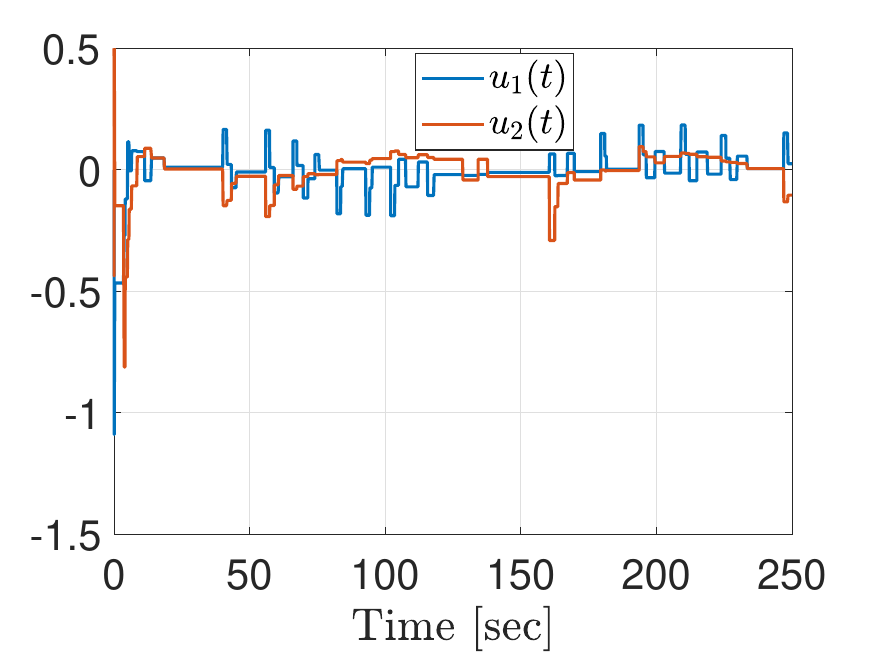}}
	\begin{picture}(0,0)
		\put(35,37){\includegraphics[width=2.2cm]{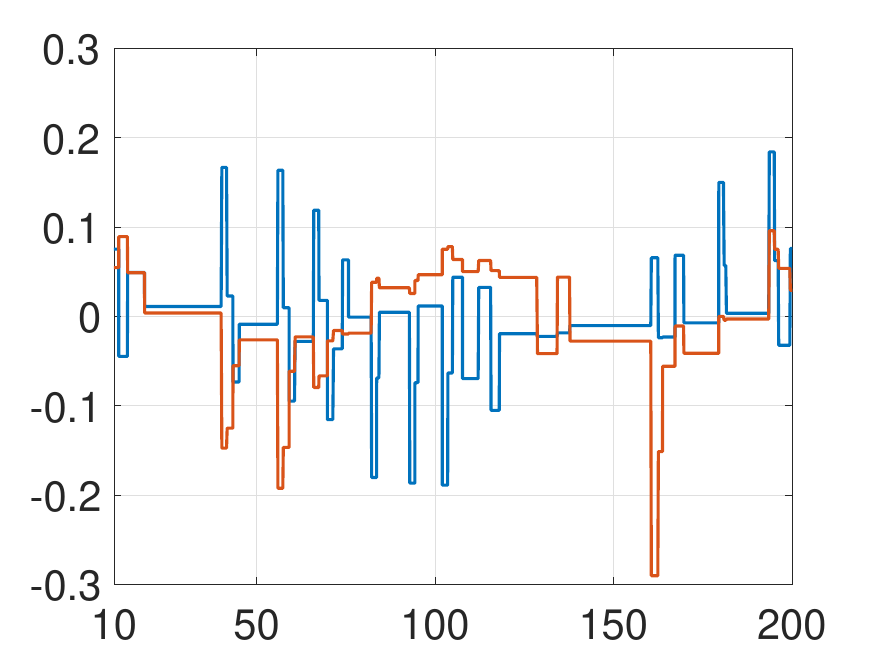}}
		\put(-75,160){\includegraphics[width=2.15cm]{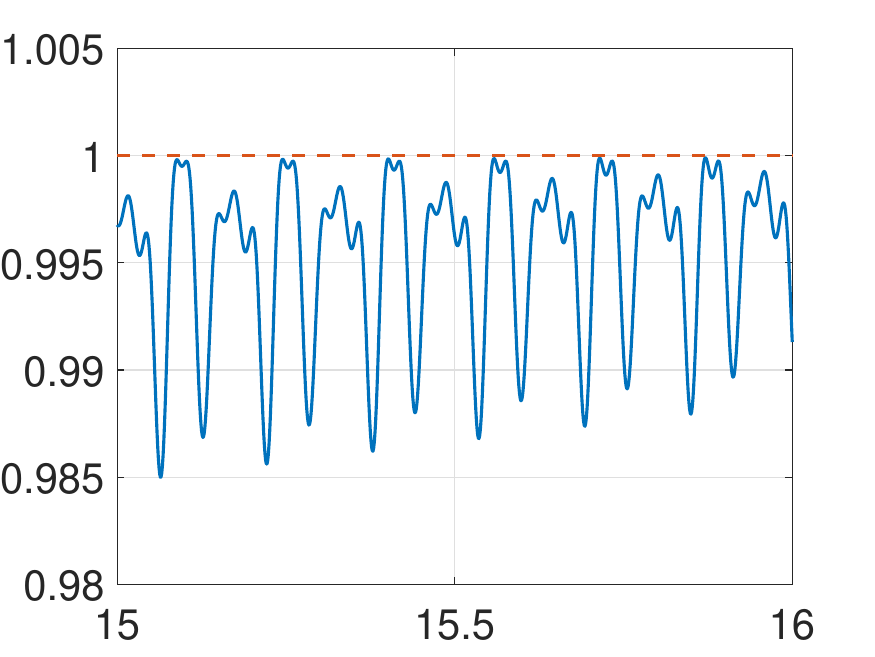}}
	\end{picture}
	\caption{Event-Triggered Tuning Law: mobile source.}
	\label{fig:ET_SSC_M}
\end{figure}

We first consider a stationary source located at $[x^{*}, y^{*}] = [2, 2]$ meters, with the vehicle starting at the origin with initial heading $\theta(0) = 0$ rad. We then let the source trace the figure-eight (Lissajous) trajectory used in \cite{SCL_2007}, namely $x^{*}(t) = a_m \sin(\omega_m t)$ and $y^{*}(t) = a_m \sin(2\omega_m t)$, with $a_m = 1$ and $\omega_m = 0.03$ rad/s $\ll \omega_0$, with the vehicle starting at $[x(0), y(0)] = [2, 2]$ meters. Fig.~\ref{fig:SSC} shows the responses of the continuous-time tuning law \eqref{eq:u}, and Figs.~\ref{fig:ET_SSC_S} and~\ref{fig:ET_SSC_M} those of the proposed ET-SSC \eqref{eq:u_ET_v1}.

\textit{Continuous tuning law (Fig.~\ref{fig:SSC}):} For the stationary source, the vehicle executes the star-shaped exploratory motion induced by the forward-velocity dither and converges to a small closed pattern around $[2,2]$ (Fig.~\ref{fig:SSC_S_xy}), while $u_1(t)$ and $u_2(t)$ peak near $1.5$ and $1.2$ and decay to zero by $t \approx 7$ s (Fig.~\ref{fig:SSC_S_u}). For the mobile source, the vehicle tracks the figure-eight path (Fig.~\ref{fig:SSC_M_xy}) with $u_1(t)$ and $u_2(t)$ settling into a small persistent oscillation, within $\pm0.1$ (Fig.~\ref{fig:SSC_M_u}).
 
\textit{Event-triggered tuning law (Figs.~\ref{fig:ET_SSC_S} and \ref{fig:ET_SSC_M}):} For both sources, the output $f(x,y)$ (Figs.~\ref{fig:ET_SSC_S_f} and \ref{fig:ET_SSC_M_f}) and the trajectories (Figs.~\ref{fig:ET_SSC_S_xy} and \ref{fig:ET_SSC_M_xy}) closely reproduce the continuous-time behavior, converging to a neighborhood of $f^{*}$ within about $5$ s; for the mobile case, Fig.~\ref{fig:ET_SSC_M_p} confirms the tracking of $x^{*}(t)$ and $y^{*}(t)$ with a small lag. The key distinction lies in the tuning signals (Figs.~\ref{fig:ET_SSC_S_u} and \ref{fig:ET_SSC_M_u}), which are held piecewise constant by the ZOH between triggering instants -- e.g., $u_1(t)$ and $u_2(t)$ are held near $0.45$ and $0.63$ from $t \approx 0.6$ s to $t \approx 2.9$ s -- in contrast to the continuously varying signals of Figs.~\ref{fig:SSC_S_u} and \ref{fig:SSC_M_u}.

\textit{Inter-event times:} In the fixed-source simulation of Fig.~\ref{fig:ET_SSC_S}, the tuning signal is updated $71$ times in $10$ s, whereas the continuous-time law is updated at each of the $10^{4}$ integration steps; in the mobile-source simulation of Fig.~\ref{fig:ET_SSC_M}, it is updated $70$ times in $250$ s, i.e., on average once every $3.6$ s, which corresponds to about eleven probing periods $T\approx 0.314$ s. For the fixed source, the inter-event times range from $2\times 10^{-3}$ s, i.e., two integration steps, which occurs only within the first $0.35$ s, when the filtered gradient estimate varies fastest, up to $4.2$ s; for the mobile source, they range from $7\times 10^{-3}$ s up to $22.8$ s. For comparison, the guaranteed lower bound \eqref{eq:F_dwell} is $\tau^{\ast}=\varepsilon_{0}/(\alpha L_{D})\approx 6.3\times 10^{-4}$ s, with $L_{D}\approx 78.8$, which is a conservative worst-case bound.


\section{CONCLUSION}

This paper developed an event-triggered source seeking control (ET-SSC)
framework for nonholonomic vehicles operating without position
measurements or prior knowledge of the source location. The proposed
design exploits a structural separation between the fast, continuous
probing motion of the vehicle and the event-triggered updating of the
controller-side source-seeking correction.
Rather than sampling the probing mechanism itself, the event-triggering
strategy acts only on this correction, allowing its computation and
transmission to occur aperiodically while the information-generating
motion remains uninterrupted.

The gradient estimate is produced by a high-pass (washout) filter
followed by a demodulation stage, and the static triggering condition,
which combines a relative and an absolute threshold, acts on its
low-pass-filtered version. Stability was proved with a
time-periodic Lyapunov function constructed along the exponentially
stable periodic solution of the filtered continuous-time closed loop,
which is continuous across the events, and practical exponential
convergence of the vehicle trajectory to a neighborhood of the unknown
source was established for both the continuous-time and the
event-triggered tuning laws. A strictly positive minimum dwell time
between consecutive events, independent of the probing frequency, was
also derived, thereby excluding Zeno behavior. The controller analyzed
and the controller simulated are identical.

Beyond the event-triggered design itself, the results highlight a distinction between designs for nonholonomic source seeking.
Steering-based and Lie-bracket source seekers offer attractive vehicle
motions and powerful analytical properties, but their gradient-producing
mechanism is intrinsically coupled to the fast motion of the vehicle.
The forward-velocity design considered here instead permits the
controller-side correction to be held between events without
interrupting the continuous probing signals. This feature makes it
particularly suitable when computation or communication of the
source-seeking correction constitutes the constrained resource.
Moreover, the averaging-based analysis employed here accommodates
anisotropic quadratic source fields with distinct spatial curvatures,
since the Hessian of the source field enters the average closed loop
directly, complementing existing Lie-bracket analyses of this design.


\section*{REFERENCES}

\def\refname{}

\vspace{-0.1cm}
\begin{IEEEbiography}[{\includegraphics[width=1in,height=1.25in,clip,keepaspectratio]{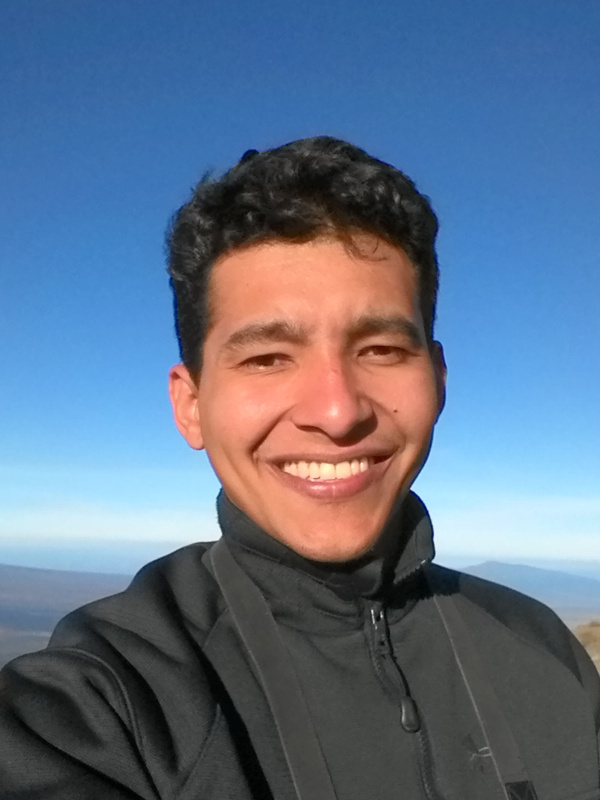}}]{VICTOR HUGO PEREIRA RODRIGUES}{\space} received the B.Sc. and M.Sc. degrees in Electronics Engineering from the State University of Rio de Janeiro (UERJ) in 2016 and 2018, respectively, and the D.Sc. degree in Electrical Engineering from the Federal University of Rio de Janeiro in 2022. In 2020, he received the FAPERJ Bolsa Nota 10 (Highest Rank Scholarship Prize). Since 2024, he has been an Associate Professor at UERJ.
\end{IEEEbiography}

\begin{IEEEbiography}[{\includegraphics[width=1in,height=1.25in,clip,keepaspectratio]{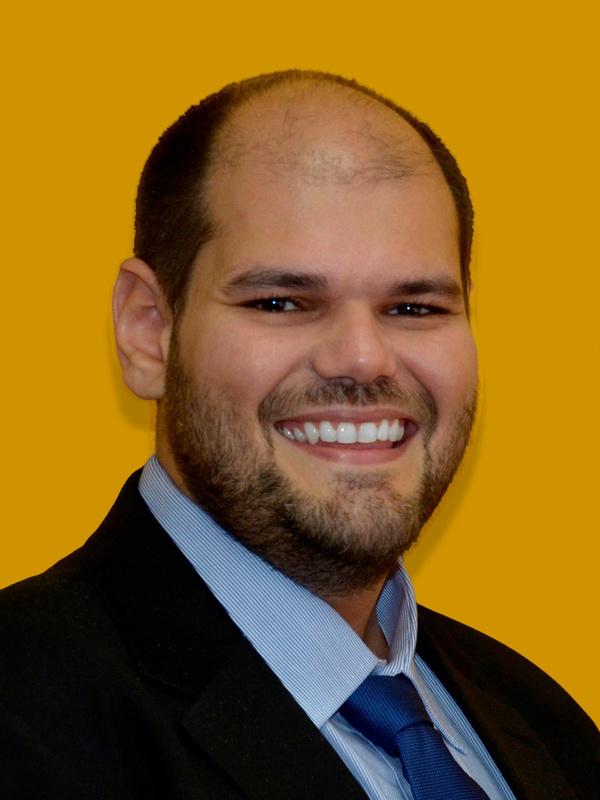}}]{TIAGO ROUX OLIVEIRA}{\space} (Senior Member, IEEE) is a Professor at UERJ and the author of over 300 publications, including the SIAM book Extremum Seeking through Delays and PDEs (2022). He chaired the IFAC TC on Adaptive and Learning Systems (2020–2026) and was President of the Brazilian Society of Automatics (2023–2025). He received the 2021 IEEE TCST Outstanding Paper Award.
\end{IEEEbiography}

\begin{IEEEbiography}[{\includegraphics[width=1in,height=1.25in,clip,keepaspectratio]{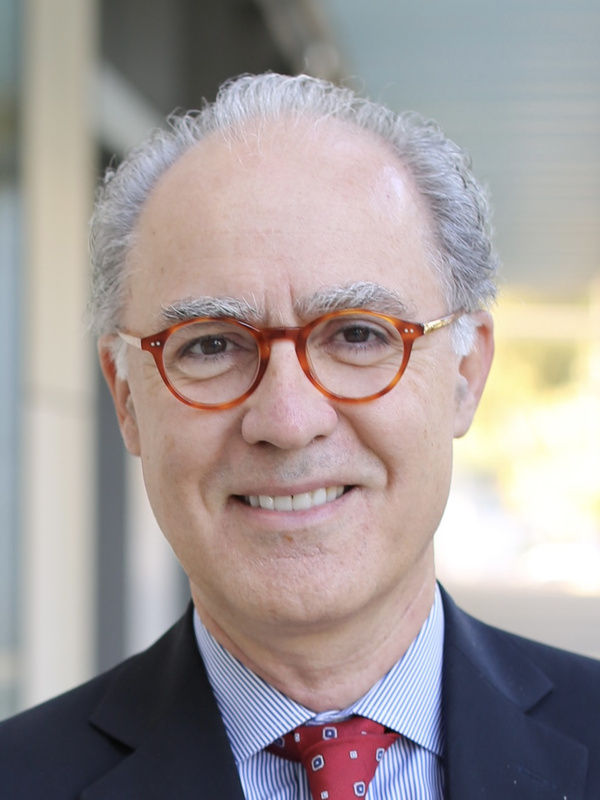}}]{MIROSLAV KRSTI{\' C}}{\space} is a Fellow of SIAM, IEEE, IFAC, ASME, AAAS, IET, AIAA (AF), and the Serbian Academy of Sciences and Arts. His awards include the SIAM Reid Prize, Bellman, Oldenburger, Ragazzini, Nyquist Lecture, Bode Lecture, and three IFAC awards. He has held chief or senior editorial positions for IEEE TAC, SCL, and Automatica, and has co-authored 19 books in control theory.
\end{IEEEbiography}

\end{document}